\documentclass{amsart}
\usepackage{amssymb}
\usepackage{hyperref}
\usepackage{mathrsfs}

\usepackage{tikz}
\usetikzlibrary{cd}

\usepackage{amssymb}
\usepackage{mathrsfs}
\usepackage{tikz-cd}

\usepackage{amsopn}

\newcommand{\G}{\mathbf{G}}
\newcommand{\B}{\mathbf{B}}
\newcommand{\T}{\mathbf{T}}
\newcommand{\F}{\mathbb{F}}
\newcommand{\X}{\mathbb{X}}
\newcommand{\Y}{\mathbb{Y}}
\newcommand{\Z}{\mathbb{Z}}
\newcommand{\Rep}{\mathsf{Rep}}
\newcommand{\Tilt}{\mathsf{Tilt}}
\newcommand{\simto}{\xrightarrow{\sim}}
\newcommand{\weyl}{\mathsf{M}}
\newcommand{\coweyl}{\mathsf{N}}
\newcommand{\simp}{\mathsf{L}}
\newcommand{\til}{\mathsf{T}}
\newcommand{\Ind}{\mathrm{Ind}}
\newcommand{\Ext}{\mathrm{Ext}}
\newcommand{\ch}{\mathrm{ch}}
\newcommand{\Waff}{W_{\mathrm{aff}}}
\newcommand{\Haff}{\mathcal{H}_{\mathrm{aff}}}
\newcommand{\Saff}{S_{\mathrm{aff}}}
\newcommand{\Wext}{W_{\mathrm{ext}}}
\newcommand{\uH}{\underline{H}}
\newcommand{\Spr}{\widetilde{\mathcal{N}}}
\newcommand{\Coh}{\mathsf{Coh}}
\newcommand{\Perv}{\mathsf{Perv}}
\newcommand{\Db}{D^{\mathrm{b}}}
\DeclareMathOperator{\Hom}{Hom}
\newcommand{\Loop}{\mathrm{L}}
\newcommand{\Gr}{\mathsf{Gr}}
\newcommand{\Iwu}{\mathrm{I}_{\mathrm{u}}}
\newcommand{\Whit}{\mathrm{Whit}}
\newcommand{\IC}{\mathscr{I} \hspace{-1pt} \mathscr{C}}
\newcommand{\U}{\mathsf{U}}
\newcommand{\C}{\mathbb{C}}
\newcommand{\su}{\mathsf{u}}

\newcommand{\Ga}{\mathbb{G}_{\mathrm{a}}}
\newcommand{\Gm}{\mathbb{G}_{\mathrm{m}}}
\newcommand{\fg}{\mathfrak{g}}
\newcommand{\cN}{\mathcal{N}}
\newcommand{\coh}{\mathsf{H}}
\newcommand{\Fr}{\mathrm{Fr}}
\DeclareMathOperator{\Spec}{Spec}
\DeclareMathOperator{\supp}{supp}
\newcommand{\oV}{\overline{V}}
\newcommand{\LR}{{\mathrm{LR}}}
\newcommand{\red}{{\mathrm{red}}}
\newcommand{\unip}{{\mathrm{unip}}}
\newcommand{\fS}{\mathfrak{S}}
\newcommand{\cL}{\mathcal{L}}
\newcommand{\cS}{\mathcal{S}}
\newcommand{\cT}{\mathcal{T}}

\numberwithin{equation}{section}
\newtheorem{thm}{Theorem}[section]

\newtheorem{conj}[thm]{Conjecture}

\theoremstyle{definition}

\theoremstyle{remark}
\newtheorem{remark}[thm]{Remark}

\title[Tilting modules for reductive algebraic groups]{Tilting modules for reductive algebraic groups: characters and support varieties}
\author{Pramod N. Achar}
\address{Department of Mathematics\\
  Louisiana State University\\
  Baton Rouge, LA 70803\\
  U.S.A.}
\email{pramod@math.lsu.edu}

\author{Simon Riche}
\address{Universit\'e Clermont Auvergne, CNRS, LMBP, F-63000 Clermont-Ferrand, France.}
\email{simon.riche@uca.fr}

\begin{document}

\maketitle

% Copyright Statement
% When submitting your final paper to a SIAM proceedings, it is requested that you include
% the appropriate copyright in the footer of the paper.  The copyright added should be
% consistent with the copyright selected on the copyright form submitted with the paper.
% Please note that "20XX" should be changed to the year of the meeting.

% Default Copyright Statement
%\fancyfoot[R]{\scriptsize{Copyright \textcopyright\ 20XX by SIAM\\
%Unauthorized reproduction of this article is prohibited}}

% Depending on which copyright you agree to when you sign the copyright form, the copyright
% can be changed to one of the following after commenting out the default copyright statement
% above.

%\fancyfoot[R]{\scriptsize{Copyright \textcopyright\ 20XX\\
%Copyright for this paper is retained by authors}}

%\fancyfoot[R]{\scriptsize{Copyright \textcopyright\ 20XX\\
%Copyright retained by principal author's organization}}

%\pagenumbering{arabic}
%\setcounter{page}{1}%Leave this line commented out.

\begin{abstract} 
These notes are our contribution to the Proceedings of the ICM 2026. We discuss some results we have obtained (in part jointly with coauthors) regarding the representation theory of reductive algebraic groups over algebraically closed fields of positive characteristic. These statements mainly concern tilting modules, in particular their characters and support varieties.
\end{abstract}

%%%%%%%%%%%%%%%%%%%%%%%%%%%%%%%%
\section{Introduction}
%%%%%%%%%%%%%%%%%%%%%%%%%%%%%%%%

%-----------------------------------------------------------------------
\subsection{Representations of reductive algebraic groups and Kazhdan--Lusztig combinatorics of the affine Weyl group}
%-----------------------------------------------------------------------

It was suggested in the early 1970s by Verma~\cite{verma} that the combinatorial invariants one can hope to compute from the representation theory of a connected reductive algebraic group $\G$ defined over an algebraically closed field of positive characteristic $p$ are closely tied to the combinatorics and geometry of the associated affine Weyl group $\Waff$. This idea was extremely influential; it was realized in various ways in work of Jantzen and Andersen, and in Lusztig's conjecture about characters of simple modules~\cite{lusztig-pbs}. Here the formula proposed by Lusztig involves the Kazhdan--Lusztig polynomials attached to $\Waff$, and it is now known to be true in large characteristics; see~\S\ref{ss:characters-simples} for some comments, and references. 

In this document we present further results in this direction obtained in the last few years, partly in collaboration with various coauthors, which exploit the idea that in fact the Kazhdan--Lusztig combinatorics only gives an ``asymptotic'' approximation of the combinatorics involved, valid for large characteristics, but that to capture the behaviour of these objects in the maximal generality one needs to consider the $p$-canonical combinatorics instead. This idea was suggested by Williamson, who had introduced the $p$-canonical basis a few years earlier in joint work with Juteau and Mautner, after his discovery of counterexamples to Lusztig's conjecture in its expected range of validity, namely $p \geq h$ with $h$ the Coxeter number of $\G$; see~\cite{williamson-torsion}. Its realization into theorems in the following years was made possible by earlier works of Bezrukavnikov, who had developed (with several coauthors) a geometric approach to the study of representations of quantum groups at roots of unity that we largely adapted to the setting of reductive groups in positive characteristic.

%-----------------------------------------------------------------------
\subsection{Characters}
%-----------------------------------------------------------------------

The first group of papers that exploited this idea, discussed in Section~\ref{sec:characters}, was concerned with the question of computing characters of appropriate representations. The character of a representation records the action of the maximal torus on this representation. It is the analogue in this setting of the character of a complex representation of a finite group, and shares many of its properties; in particular it determines the composition factors of the representation (and their multiplicities), provided we know the characters of simple representations. It is therefore a central question to compute the latter characters. More specifically, there are nice representations for which characters are known, namely induced modules, and one expects to express characters of simple modules as linear combinations of characters of these modules; this is exactly the setting of Lusztig's conjecture.

Our approach to this question follows work of Andersen in the early 1990s, who remarked that (under mild assumptions on $p$) an expression for characters of simple modules can be derived from an expression for characters of (some) modules in another family parametrized in a similar way, that of indecomposable \emph{tilting} representations. Andersen also proposed a conjectural answer to this question for a suitable subfamily, and showed that it implied Lusztig's conjecture. Williamson and the second author conjectured in~\cite{rw-tilting} that, replacing in Andersen's conjecture the Kazhdan--Lusztig polynomials by their $p$-versions (i.e.~the polynomials obtained by the same procedure using the $p$-canonical basis instead of the Kazhdan--Lusztig basis), one obtains a formula for characters of \emph{all} tilting modules, valid as soon as $p>h$. This conjecture was first proved (based on adaptations of several constructions of Bezrukavnikov) by Makisumi, Williamson and the authors in~\cite{amrw}; it now admits several independent proofs, including a version (obtained by Williamson and the second author) that makes sense (and holds) for all values of $p$. See~\S\ref{ss:tilting-char} for a discussion of these results, and~\cite{williamson-alg,williamson} for other presentations.

As explained above, from the knowledge of characters of indecomposable tilting modules one can in theory derive characters of simple modules, but this procedure does not directly leads to explicit formulas. Further efforts in this direction are therefore useful. Here we explain an answer to this question (when $p>h+1$) that can be derived from a conjecture due to Finkelberg--Mirkovi\'c, proved recently by Bezrukavnikov and the second author. The formula so obtained involves some geometric data that are currently not well understood, namely Euler characteristics of stalks of Iwahori-equivariant intersection cohomology complexes on the affine flag variety of the Langlands dual group; see~\S\ref{ss:FM-conj} for details.

%-----------------------------------------------------------------------
\subsection{Support varieties}
%-----------------------------------------------------------------------

The second group of publications that we discuss (in Section~\ref{sec:support}) is concerned with the description of support varieties of indecomposable tilting modules. This question is the subject of a celebrated conjecture due to Humphreys~\cite{humphreys}, that we have proved in some cases in collaboration with Hardesty. The answer proposed by Humphreys also involves the Kazhdan--Lusztig combinatorics of $\Waff$, and more specifically its two-sided cells, which (thanks to work of Lusztig~\cite{lusztig-cawg4}) are in a canonical bijection with nilpotent orbits for (the Frobenius twist of) $\G$. The analogous question for quantum groups at a root of unity was studied by Ostrik, and solved in full generality by Bezrukavnikov~\cite{bezrukavnikov}. In fact Bezrukavnikov proved in this case a more precise statement, relating (relative) cohomology of indecomposable tilting modules to perverse-coherent sheaves on the nilpotent cone of the corresponding complex reductive group, via a refinement of Lusztig's bijection called the Lusztig--Vogan bijection.

Since the main tools used in~\cite{bezrukavnikov} had been adapted to the setting of reductive groups in positive characteristic in the course of our work on characters of tilting modules, we worked (with Hardesty) towards an adaptation of this approach in~\cite{ahr}. We however had to face two main difficulties. The first one is that, as explained above, the combinatorics that appears in our work is that of the $p$-canonical basis rather than the Kazhdan--Lusztig combinatorics; the cells that naturally occur in this study are therefore the $p$-cells (the analogues of Kazhdan--Lusztig cells for the $p$-canonical basis). In this way the Humphreys conjecture translates into subtle compatibility properties of $p$-cells with Kazhdan--Lusztig cells, which we can establish for now only under the assumption that $p$ is large (with no explicit bound), leading to a proof of the Humphreys conjecture in this setting. See~\cite{ahr-conj} for a more detailed presentation of these results, and a conjectural broader picture including them.

The second difficulty is that, for reductive groups in positive characteristic, the coherent sheaves on the nilpotent cone arising from relative cohomology of tilting modules are \emph{not} the simple perverse-coherent sheaves. A solution to this problem was found by Hardesty and the first author in~\cite{ah1,ah2}, where it was proved that the appropriate framework for the description of these objects is that of co-t-structures. More precisely, in~\cite{ah1} it is proved that the objects under consideration are the indecomposable objects in the coheart of a certain co-t-structure on the derived category of equivariant coherent sheaves on the nilpotent cone of (the Frobenius twist of) $\G$, and in~\cite{ah2} a different construction of this co-t-structure is given, making it possible to describe the support of these objects in terms of the Lusztig--Vogan bijection, and leading to a proof of a relative variant of Humphreys' conjecture under the assumption that $p>h$.

%%%%%%%%%%%%%%%%%%%%%%%%%%%%%%%%
\section{Characters}
\label{sec:characters}
%%%%%%%%%%%%%%%%%%%%%%%%%%%%%%%%

%-----------------------------------------------------------------------
\subsection{Representations of reductive algebraic groups and their characters}
%-----------------------------------------------------------------------

Let $\Bbbk$ be an algebraically closed field of characteristic $p$, and let $\G$ be a connected reductive algebraic group over $\Bbbk$. We also fix a Borel subgroup $\B \subset \G$ and a maximal torus $\T \subset \B$, and denote by $\X$, resp.~$\Y$, the lattice of weights, resp.~coweights, of $\T$, i.e.~morphisms of $\Bbbk$-algebraic groups from $\T$ to the multiplicative group $\Gm$, resp.~from $\Gm$ to $\T$. Then $\X$ and $\Y$ are naturally dual lattices, and in $\X$, resp.~$\Y$, we have the roots $\Phi$, resp.~coroots $\Phi^\vee$. The choice of $\B$ determines subsets $\Phi_+ \subset \Phi$ and $\Phi^\vee_+ \subset \Phi^\vee$ of positive roots and coroots, where $\Phi_+$ consists of the roots that are opposite to the $\T$-weights in the Lie algebra of $\B$. One can then define the subset $\X_+ \subset \X$ of dominant weights to consist of the weights whose pairing with any element of $\Phi^\vee_+$ is nonnegative. We will also consider the natural order $\preceq$ on $\X$, such that $\lambda \preceq \mu$ iff $\mu-\lambda$ is a sum of elements of $\Phi_+$. Finally, we will denote by $W$ the Weyl group of $\G$ (with respect to $\T$); this group naturally acts on $\T$, hence also on $\X$. The choice of $\B$ (or of $\Phi_+$) determines a subset $S \subset W$ of Coxeter generators, and hence in particular a length function $\ell$ on $W$.

Our main object of study in this paper will be the category $\Rep(\G)$ of finite-dimensional algebraic representations of $\G$, whose Grothendieck group will be denoted $[\Rep(\G)]$. Tensor product of representations defines a monoidal structure on $\Rep(\G)$, and thus also a ring structure on $[\Rep(\G)]$. The most interesting combinatorial invariant one can attach to an object of $\Rep(\G)$ is its \emph{character}, defined as follows. From the representation theory of tori we know that for any $M$ in $\Rep(\G)$ we have
\[
M = \bigoplus_{\lambda \in \X} M_\lambda \quad \text{where $M_\lambda = \{v \in M \mid \forall t \in \T, \, t \cdot v = \lambda(t) v\}$;}
\]
then we denote by $\ch(M)$ the formal sum $\sum_{\lambda \in \X} \dim(M_\lambda) \cdot e^\lambda \in \Z[\X]$. It is well known that the assignment $M \mapsto \ch(M)$ induces a ring isomorphism $[\Rep(\G)] \simto \Z[\X]^W$.

\begin{remark}
\phantomsection
\label{rmk:gps}
\begin{enumerate}
\item
All the questions we will consider in this section can easily be reduced to the case where $\G$ has simply connected derived subgroup.
%; this assumption is sometimes imposed in the references that we give, but we will not impose it here if it is not necessary.
\item
\label{it:gps-fields}
It is known that connected reductive algebraic groups come ``in families over all fields.'' More formally, for any group $\G$ as above there exists a unique (up to isomorphism) split reductive group scheme $\G_{\Z}$ over $\Z$ such that $\G = \Spec(\Bbbk) \times_{\Spec(\Z)} \G_{\Z}$. One can also choose a split maximal torus and a Borel subgroup $\T_\Z \subset \B_\Z \subset \G_\Z$ and choose the subgroups $\T$ and $\B$ as above as $\B = \Spec(\Bbbk) \times_{\Spec(\Z)} \B_{\Z}$, $\T = \Spec(\Bbbk) \times_{\Spec(\Z)} \T_{\Z}$. Then, for any algebraically closed field $\mathbb{L}$ one can consider the connected reductive algebraic group $\G_{\mathbb{L}} = \Spec(\mathbb{L}) \times_{\Spec(\Z)} \G_{\Z}$, its Borel subgroup $\B_{\mathbb{L}} = \Spec(\mathbb{L}) \times_{\Spec(\Z)} \B_{\Z}$, its maximal torus $\T_{\mathbb{L}} = \Spec(\mathbb{L}) \times_{\Spec(\Z)} \T_{\Z}$, and ask to what extent the category $\Rep(\G_{\mathbb{L}})$ depends on the choice of field $\mathbb{L}$. (For a general group $\G$, the construction of $\G_\Z$ relies on a lot of difficult theory, but for many groups this is something well known: it is clear that general linear groups, symplectic groups or special orthogonal groups can be considered over any field.) We will come back to this point of view in various remarks.
\end{enumerate}
\end{remark}

%-----------------------------------------------------------------------
\subsection{Induced, simple and tilting modules}
\label{ss:ind-simp-tilt}
%-----------------------------------------------------------------------

The category $\Rep(\G)$ contains a number of families of objects of particular interest, that we now introduce. It is a standard fact that any $\lambda \in \X$ extends in a unique way to a morphism from $\B$ to $\Gm$; to $\lambda$ one can therefore associate a $1$-dimensional $\B$-module $\Bbbk_{\B}(\lambda)$, that we can then induce to $\G$ to obtain the algebraic $\G$-module $\coweyl(\lambda) = \Ind_{\B}^{\G}(\Bbbk_{\B}(\lambda))$. The following properties of these modules, which we will call the induced modules, are well known. (All the results described here and below as ``well known'' are discussed e.g.~in the seminal book of Jantzen~\cite{jantzen}.)
\begin{enumerate}
\item
For any $\lambda \in \X$, $\coweyl(\lambda)$ is finite-dimensional, i.e.~an object of $\Rep(\G)$.
Moreover, $\coweyl(\lambda) \neq 0$ iff $\lambda \in \X_+$.
\item
If $\lambda \in \X_+$, $\coweyl(\lambda)$ contains a unique simple sub-$\G$-module $\simp(\lambda) \subset \coweyl(\lambda)$.
Moreover,
the assignment $\lambda \mapsto \simp(\lambda)$ induces a bijection between $\X_+$ and the set of isomorphism classes of simple objects in $\Rep(\G)$.
\end{enumerate}
For $\lambda \in \X_+$ we will also denote by $\weyl(\lambda)$ the $\G$-module dual to $\coweyl(-w_\circ \lambda)$ (where $w_\circ$ is the longest element in $W$); these modules are called the Weyl modules.

The last family of modules that we will consider, and that will be the main focus of this text, is that of indecomposable tilting modules. An object $M$ of $\Rep(\G)$ is called \emph{tilting} if both $M$ and the dual module $M^*$ admit (finite) filtrations whose subquotients are of the form $\coweyl(\mu)$ for some $\mu \in \X_+$. It is clear that this class of modules is stable under direct sums, and it is a classical fact that it is also stable under direct summands, which reduces their description to that of the indecomposable tilting modules. Such modules are again classified by $\X_+$; more specifically:
\begin{enumerate}
\item
 for any $\lambda \in \X_+$ there exists a unique indecomposable tilting module $\til(\lambda)$ such that $\dim(\til(\lambda)_\lambda)=1$ and $\til(\lambda)_\mu=0$ unless $\mu \preceq \lambda$;
 \item
  the assignment $\lambda \mapsto \til(\lambda)$ induces a bijection between $\X_+$ and the set of isomorphism classes of indecomposable tilting $\G$-modules. 
 \end{enumerate}
It is a fundamental (and quite nontrivial) fact that a tensor product of tilting $\G$-modules is again tilting; the full subcategory $\Tilt(\G) \subset \Rep(\G)$ consisting of tilting modules is therefore a monoidal subcategory.

In case $p=0$, it turns out that $\simp(\lambda)=\coweyl(\lambda)=\weyl(\lambda)=\til(\lambda)$ for any $\lambda \in \X_+$, but when $p>0$ these classes \emph{are} different (unless $\G$ is a torus). The simple modules are the ``building blocks'' of all representations, in the sense that any object of $\Rep(\G)$ is an extension of such modules; their study is thus of particular importance. Although their definition might seem convoluted at first, the class of tilting modules turns out to be also extremely interesting; in particular, due to results of Sobaje~\cite{sobaje} (refining earlier work of Andersen~\cite{andersen-tilting} which imposed a technical condition on $p$), there is a procedure to describe the characters of simple modules if the characters of some indecomposable tilting modules are known. (This procedure is explicit but not really straightforward, and will not be discussed here in detail.)

%-----------------------------------------------------------------------
\subsection{The Weyl character formula}
\label{ss:Weyl-formula}
%-----------------------------------------------------------------------

For each of the three families of modules introduced in~\S\ref{ss:ind-simp-tilt}, one can ask the question of describing the characters of its constituents. The answer to this question is well understood for induced modules and Weyl modules, and given by a formula identical to one due to Weyl for compact Lie groups: for any $\lambda \in \X_+$ we have
\begin{equation}
\label{eqn:Weyl-formula}
\ch(\coweyl(\lambda)) = \ch(\weyl(\lambda)) = \frac{\sum_{w \in W} (-1)^{\ell(w)} \cdot e^{w(\lambda + \rho)-\rho}}{\sum_{w \in W} (-1)^{\ell(w)} \cdot e^{w(\rho)-\rho}}
\end{equation}
where $\rho$ is one-half the sum of the positive roots. (In this formula, all exponents belong to $\X$, although $\rho$ is only an element of $\frac{1}{2}\X$ in general.)

Given this information, and since the character of a module only depends on its class in $[\Rep(\G)]$, for any $M$ in $\Rep(\G)$, to describe $\ch(M)$ it suffices to express its class $[M]$ as a linear combination of the classes of the induced modules. (It turns out that $([\coweyl(\lambda)] : \lambda \in \X_+)$ is a basis of $[\Rep(\G)]$, so that such an expansion always uniquely exists.) This is the reformulation that we will consider below for the cases of simple and indecomposable tilting modules.

\begin{remark}
\begin{enumerate}
\item
From the definitions, it is clear that the coefficients of the expansion of each  $\til(\mu)$ in the basis $( [\coweyl(\lambda)] : \lambda \in \X_+ )$ are nonnegative. For simple modules however, some coefficients might be negative.
\item
Recall the point of view introduced in Remark~\ref{rmk:gps}\eqref{it:gps-fields}. A first (naive) attempt at using this perspective is as follows. Using the notation of that remark, for any algebraically closed field $\mathbb{L}$ the character lattice of $\T_{\mathbb{L}}$ identifies canonically with $\X$, and the subset of dominant weights $\X_+$ is also the same for all fields. Given $\lambda \in \X_+$, one can consider the representations $\simp_{\mathbb{L}}(\lambda)$, $\coweyl_{\mathbb{L}}(\lambda)$ or $\til_{\mathbb{L}}(\lambda)$ over $\mathbb{L}$, and ask if their characters vary with $\mathbb{L}$. It is easily seen that these characters only depend on the characteristic of $\mathbb{L}$, and for the case of induced modules $\coweyl(\lambda)$, the Weyl character formula~\eqref{eqn:Weyl-formula} shows that the answer does not depend on $\mathbb{L}$ at all. It can also be shown that if we fix $\lambda$, if the characteristic of $\mathbb{L}$ is larger than a certain bound depending on $\lambda$, we have $\simp_{\mathbb{L}}(\lambda)=\coweyl_{\mathbb{L}}(\lambda)=\weyl_{\mathbb{L}}(\lambda)=\til_{\mathbb{L}}(\lambda)$.
\end{enumerate}
\end{remark}

%-----------------------------------------------------------------------
\subsection{Linkage principle}
\label{ss:linkage}
%-----------------------------------------------------------------------

From now on we will assume that $p>0$. In this case,
as was first remarked by Verma~\cite{verma},
the structure of the category $\Rep(\G)$ exhibits some kinds of symmetries governed by the affine Weyl group $\Waff := W \ltimes \Z\Phi$ (where $\Z\Phi \subset \X$ is the sublattice generated by $\Phi$) and its action on $\X$ (called the ``dot-action'') defined by
\begin{equation}
\label{eqn:dot-action}
(w t_\mu) \bullet \lambda = w(\lambda+p\mu+\rho)-\rho
\end{equation}
for $w \in W$, $\mu \in \Z\Phi$ and $\lambda \in \X$. (Here, we use the notation $t_\mu$ for the element $(1,\mu) \in \Waff$.) More specifically, the \emph{linkage principle}, conjectured by Verma and proved in full generality by Andersen~\cite{andersen} (after earlier proofs under additional technical assumptions by Humphreys, Jantzen and Carter--Lusztig) states that for $\lambda, \mu \in \X_+$ we have
\[
\Ext^1_{\Rep(G)}(\simp(\lambda),\simp(\mu))=0 \quad \text{unless $\Waff \bullet \lambda = \Waff \bullet \mu$.}
\]
As a consequence, if for any subset $c \subset \X$ we denote by $\Rep_c(\G) \subset \Rep(\G)$ the full subcategory whose objects are the modules whose composition factors all belong to $\{ \simp(\lambda) : \lambda \in c \}$, we have
\begin{equation}
\label{eqn:Rep-blocks}
\Rep(\G) = \bigoplus_{c \in \X/(\Waff,\bullet)} \Rep_c(\G).
\end{equation}

\begin{remark}
\begin{enumerate}
\item
For $c \in \X/(\Waff,\bullet)$,
the subcategory $\Rep_c(\G)$ is often called the ``block of $c$.'' One should be careful that this subcategory is not a ``block'' in the usual sense, because it might be decomposable.
\item
The interplay between the monoidal structure on $\Rep(\G)$ and the decomposition~\eqref{eqn:Rep-blocks} is a very interesting subject, which however is not well understood. In particular, the tensor product of two modules which each belong to one block usually does not belong to one block.
\item
The affine Weyl group $\Waff$ is a subgroup of the ``extended affine Weyl group'' $\Wext := W \ltimes \X$. The formula~\eqref{eqn:dot-action} makes sense more generally for any $\mu \in \X$, which allows us to extend the dot-action of $\Waff$ on $\X$ to $\Wext$; we will denote this action by the same symbol.
\item
In the point of view of Remark~\ref{rmk:gps}\eqref{it:gps-fields}, the group $\Waff$ is the same for all fields $\mathbb{L}$, but the dot-action on $\X$ \emph{does} depend on the characteristic of $\mathbb{L}$.
\end{enumerate}
\end{remark}

If $M \in \Rep(\G)$ is indecomposable, there exists a unique $c \in \X/(\Waff,\bullet)$ such that $M \in \Rep_c(\G)$. Then, in the expansion of $[M]$ as a linear combination of classes of induced modules (see~\S\ref{ss:Weyl-formula}), only the classes of modules labelled by elements of $c$ can occur. The problem of describing characters of simple or indecomposable tilting modules can therefore be reformulated as follows: describe, for any $(\Waff,\bullet)$-orbit $c \subset \X$ and any $\lambda \in \X_+$, the expansion of the classes of $\simp(\lambda)$ and $\til(\lambda)$ as a linear combination of the classes of the modules $(\coweyl(\mu) : \mu \in c \cap \X_+)$. In fact there exists a particularly nice fundamental domain for the action of $\Waff$ on $\X$, namely the subset
\[
C := \{\lambda \in \X \mid \forall \alpha^\vee \in \Phi^\vee_+, \, 0 \leq \langle \lambda+\rho, \alpha^\vee \rangle \leq p\}.
\]
The choice of an orbit is therefore equivalent to the choice of an element in $C$.

%-----------------------------------------------------------------------
\subsection{Translation principle}
\label{ss:translation}
%-----------------------------------------------------------------------

The reformulation of our problem considered at the end of~\S\ref{ss:linkage} can be simplified using further structure on the category $\Rep(\G)$, which involves a canonical Coxeter group structure on $\Waff$. Namely, there exists a canonical subset $\Saff \subset \Waff$ containing $S$ and such that the pair $(\Waff,\Saff)$ is a Coxeter system. (Here, ``canonical'' means that this subset only depends on our initial data of $\G$, $\B$, $\T$; but it does depend on the choice of $\B$.) 

\begin{remark}\label{rmk:wext-omega}
The group $\Wext$ is not naturally a Coxeter group, but it has a similar structure: the length function $\ell : \Waff \to \Z_{\geq 0}$ extends canonically to $\Wext$, in such a way that setting $\Omega = \{w \in \Wext \mid \ell(w)=0\}$, conjugation by any $\omega \in \Omega$ preserves $\Saff$, and multiplication induces a group isomorphism $\Omega \ltimes \Waff \simto \Wext$.
\end{remark}

A subset $J \subset \Saff$ is called \emph{finitary} if the subgroup $W_J$ it generates in $\Waff$ is finite; in this case we will denote by $w_J$ the longest element in $W_J$.
Given a finitary subset $J \subset \Saff$, we will denote by $C_J \subset C$ the (possibly empty) subset consisting of elements whose stablizer in $\Waff$ is $W_J$;
then we have
\[
C = \bigsqcup_{J \subset \Saff \text{ finitary}} C_J.
\]
%of $C$ into ``facets'' parametrized by finitary subsets of $\Saff$ (some of which might be empty) which 
Moreover this partition satisfies the following properties. (These statements are reformulations of results from~\cite[Chap.~II.7 \& \S II.E.11]{jantzen}; for the present point of view, see also~\cite[Chap.~1, \S 2]{riche}.)

\begin{enumerate}
\item
If $J \subset \Saff$ is finitary, 
%and $\lambda \in C_I$ one can describe the elements of $(\Waff \bullet \lambda) \cap \X_+$ in Coxeter-theoretic terms as follows. If 
we will denote by $\Waff^J$ the subset of $\Waff$ consisting of the elements $w$ which are ``strongly minimal in $W w W_J$'' in the sense that for any $x \in W$ and $y \in W_J$ we have $\ell(xwy) = \ell(x) + \ell(w) + \ell(y)$. (See~\cite[Lemma~2.4]{ar-steinberg} for other characterizations of this subset.)
%of minimal length in the coset $Ww$ and of maximal length in the coset $wW_I$. 
Then for any $\lambda \in C_J$ the assignment $w \mapsto w \bullet \lambda$ induces a bijection between $\Waff^J$ and $(\Waff \bullet \lambda) \cap \X_+$.
\item
Consider now finitary subsets $J \subset K \subset \Saff$, so that $W_J \subset W_K$. If we denote by $W_K^J \subset W_K$ the subset consisting of the elements $y$ which are minimal in $yW_J$, for any $y \in W_K^J$ we have $\Waff^K \cdot y \subset \Waff^J$. 
%Let us also denote by $w_J^I$ the unique element in $W_J^I \cap w_J W_I$. 
Given elements $\lambda \in C_J$ and $\mu \in C_K$, there exist exact biadjoint ``translation functors''
\begin{equation}
\label{eqn:translation-functors}
T_\lambda^\mu : \Rep_{\Waff \bullet \lambda}(\G) \to \Rep_{\Waff \bullet \mu}(\G), \quad T^\lambda_\mu : \Rep_{\Waff \bullet \mu}(\G) \to \Rep_{\Waff \bullet \lambda}(\G)
\end{equation}
such that for any $w \in \Waff^J$ and $y \in \Waff^K$ we have:
\begin{enumerate}
\item
\label{it:translation-coweyl}
$T_\lambda^\mu(\coweyl(w \bullet \lambda))$ is isomorphic to $\coweyl(w \bullet \mu)$ if $wW_K \cap \Waff^K \neq \varnothing$ (in other words, if $w$ is of the form $uz$ with $u \in \Waff^K$ and $z \in W_K^J$), and is $0$ otherwise;
\item
the module $T_\mu^\lambda(\coweyl(y \bullet \mu))$ admits a filtration whose associated graded is $\bigoplus_{z \in W_K^J} \coweyl(yz \bullet \lambda)$;
\item
\label{it:translation-simples}
$T_\lambda^\mu(\simp(w \bullet \lambda))$ is isomorphic to $\simp(w \bullet \mu)$ if $w \in \Waff^K$, and is $0$ otherwise;
\item
$T_\mu^\lambda(\til(y \bullet \mu)) \cong \til(y w_K w_J \bullet \lambda)$.
\end{enumerate}
In particular, when $J=K$ the functors $T_\lambda^\mu$ and $T^\lambda_\mu$ are mutually inverse equivalences which send simple, induced and indecomposable tilting modules to simple, induced and indecomposable tilting modules respectively.
\end{enumerate}

These properties have the following consequences. First, given a finitary subset $J \subset \Saff$ such that $C_J \neq \varnothing$, there exist collections of integers $(c^J_{y,w} : y,w \in \Waff^J)$ and $(d^J_{y,w} : y,w \in \Waff^J)$ such that for any $\lambda \in C_J$ and $w \in \Waff^J$ we have
\[
[\simp(w \bullet \lambda)] = \sum_{y \in \Waff^J} c^J_{y,w} \cdot [\coweyl(y \bullet \lambda)], \qquad [\til(w \bullet \lambda)] = \sum_{y \in \Waff^J} d^J_{y,w} \cdot [\coweyl(y \bullet \lambda)].
\]
(In other words, these integers, which might a priori depend on the choice of $\lambda$, in fact do not. This property, for the coefficients $c^J_{y,w}$, was conjectured by Verma, see~\cite[Conjecture~III]{verma}.) Second, given finitary subsets $J \subset K \subset \Saff$ such that $C_J \neq \varnothing$ and $C_K \neq \varnothing$, if we know the values of the coefficients $(c^J_{y,w} : y,w \in \Waff^J)$, resp.~$(d^J_{y,w} : y,w \in \Waff^J)$, then one can deduce the values of the coefficients $(c^K_{y,w} : y,w \in \Waff^K)$, resp.~$(d^K_{y,w} : y,w \in \Waff^K)$. More specifically, for $y,w \in \Waff^K$ we have
\begin{equation}
\label{eqn:formulas-c-d-translation}
c^K_{y,w} = \sum_{z \in W_K^J} c^J_{yz,w}, \qquad d^K_{y,w} = d^J_{yu,ww_K w_J} \text{ for any $u \in W_K^J$}.
\end{equation}

The most favorable situation is when $C_{\varnothing} \neq \varnothing$: in this case the formulas above show that the coefficients associated with any finitary subset $J$ are determined by those for the subset $\varnothing$ (corresponding to the so-called ``regular blocks''). Unfortunately this condition is not always satisfied; in fact we have $C_{\varnothing} \neq \varnothing$ iff $p \geq h$, where $h$ is the Coxeter number of the root system $\Phi$. In this case the weight $0$ belongs to $C_\varnothing$. 

\begin{remark}
\label{rmk:coeff-c-d}
\begin{enumerate}
\item
\label{it:tensor-product-formulas}
There is another important result regarding simple modules that we do not discuss here, called \emph{Steinberg's tensor product theorem}, which reduces, in case $\G$ has simply connected derived subgroup, the determination of their characters to the case when $\lambda$ is restricted, i.e. $\langle \lambda, \alpha^\vee \rangle < p$ for any simple coroot $\alpha^\vee$. With the parametrization above, for $\lambda \in C_J$ and $w \in \Waff^J$, $w\bullet \lambda$ is restricted iff $w^{-1}$ is restricted in the sense considered in~\cite[\S 2.4]{ar-steinberg}; in particular, only a finite number of elements $w$ need to be considered. It is not so simple, however, to express explicitly this reduction in terms of the coefficients $(c^J_{y,w} : y,w \in \Waff^J)$, due to the absence of a formula of this type for induced modules. 

There exists also a formula of a similar type for indecomposable tilting modules, due to Donkin (and known as of now only under the assumption\footnote{It was long expected that this formula should hold for all $p$. However counterexamples in small characteristics were recently found by Bendel--Nakano--Pillen--Sobaje~\cite{bnps2, bnps3}, who also improved the original bound $p \geq 2h-2$ (deduced by Donkin from earlier results of Jantzen) to $p \geq 2h-4$. The precise range of validity of this formula is still not understood.} that $p \geq 2h-4$); in this case however, it does not reduce the question to the consideration of a finite number of $w$'s.
\item
\label{it:multiplicities-indep-p}
Recall the point of view introduced in Remark~\ref{rmk:gps}\eqref{it:gps-fields}. The group $\Waff$ and its Coxeter group structure are the same for all groups $\G_{\mathbb{L}}$, and one can consider the coefficients $(c^J_{y,w} : y,w \in \Waff^J)$ and $(d^J_{y,w} : y,w \in \Waff^J)$ for all characteristics such that $C_J \neq \varnothing$. Due to the statements mentioned in~\eqref{it:tensor-product-formulas}, these coefficients must depend on the characteristic;\footnote{The independence of all coefficients $(c^J_{y,w} : y,w \in \Waff^J)$ was conjectured by Verma in~\cite[Conjecture~V]{verma}, but this expectation was too optimistic, as discussed in~\cite[\S II.8.22]{jantzen}.} but it has been a guiding principle for long that, if we restrict $w$ to live an appropriate region, this dependence should disappear when the characteristic of $\mathbb{L}$ is larger than $h$ (or at least larger than a bound linear in $h$). This idea was first suggested in~\cite{verma}; it was later made much more precise by Lusztig~\cite{lusztig-pbs}, who proposed an explicit formula for the coefficients $c_{y,w}^\varnothing$, see Conjecture~\ref{conj:lusztig} below,
%in terms of the polynomials $(h_{y,w} : y,w \in \Waff)$, and 
which is visibly independent of $p$ in its expected range of validity. 
%In fact, Lusztig's conjecture mentioned below is a precise formulation of this idea (which has been at the source of this philosophy), and 
Andersen's conjecture (also discussed below) is of the same form. It is now known that Lusztig's conjecture is true if the characteristic is larger than some bound (depending only on $\Phi$), but it is also known that this bound cannot be polynomial in $h$, see~\cite{williamson-torsion}. It has been a major motivation for our work described below to try to find formulas for these coefficients that apply for all characteristics, or at least under much milder assumptions.
\end{enumerate}
\end{remark}

%-----------------------------------------------------------------------
\subsection{Affine Hecke algebra and Kazhdan--Lusztig basis}
\label{ss:Haff}
%-----------------------------------------------------------------------

As with any Coxeter system, to $\Waff$ one can associate its Hecke algebra $\Haff$: it is a $\Z[v,v^{-1}]$-algebra with a basis $(H_w : w \in \Waff)$ and multiplication uniquely determined by the following rules:
\[
H_x \cdot H_y = H_{xy} \text{ if $\ell(xy)=\ell(x)+\ell(y)$,} \qquad (H_s+v)(H_s-v^{-1})=0 \text{ for all $s \in \Saff$.}
\]
Kazhdan--Lusztig~\cite{kazhdan-lusztig} have introduced (again, in the generality of Coxeter systems) a remarkable basis $(\uH_w : w \in \Waff)$ of this algebra, characterized by some simple algebro-combinatorial properties, which has been shown to capture important combinatorial information on the representation theory of various objects associated with $\Waff$. (Here, we use the notation and normalization of~\cite{soergel}.) This basis gives rise to the Kazhdan--Lusztig polynomials $(h_{y,w} : y,w \in \Waff)$, defined as the coefficients of the expansion of the elements $(\uH_w : w \in \Waff)$ in the basis $(H_y : y \in \Waff)$:
\[
\uH_w = \sum_{y \in \Waff} h_{y,w} \cdot H_y.
\]
(It is part of the characterization of the Kazhdan--Lusztig basis that these coefficients are polynomials in $v$.)

%\begin{remark}
%\label{rmk:multiplicities-indep-p}
%\end{remark}

More recently a new family of bases of $\Haff$ (and, more generally, of the Hecke algebra of any crystallographic Coxeter system), parametrized by prime numbers, has been introduced by Williamson:
%, see e.g.~\cite{jensen-williamson,williamson}: 
for any prime number $q$ we have the $q$-canonical basis $({}^q \hspace{-1pt} \uH_w : w \in \Waff)$, and one defines the $q$-Kazhdan--Lusztig polynomials $({}^q \hspace{-1pt} h_{y,w} : y,w \in \Waff)$ as the coefficients obtained by expanding the elements of this basis in the standard basis:
\[
{}^q \hspace{-1pt} \uH_w = \sum_{y \in \Waff} {}^q \hspace{-1pt} h_{y,w} \cdot H_y.
\]
(Here the coefficients can be Laurent polynomials rather than polynomials.)
We refer the reader to~\cite{jensen-williamson,williamson} or to~\cite[Chap.~2]{riche} for reviews of the main properties of this basis. Here we just note the following facts.
\begin{enumerate}
\item
The definition of the $q$-canonical basis involves categorification; in particular, it does not have an ``easy'' algebraic characterization just in terms of $\Haff$. More specifically one considers a certain (additive, monoidal) ``Hecke category'' $\mathscr{H}_{\mathrm{aff}}^{\mathbb{L}}$ with coefficients in a field $\mathbb{L}$ of characteristic $q$, whose split Grothendieck group identifies with $\Haff$, and defines the $q$-canonical basis as the classes of (normalized) indecomposable objects in $\mathscr{H}_{\mathrm{aff}}^{\mathbb{L}}$. (There are several possible choices for the category $\mathscr{H}_{\mathrm{aff}}^{\mathbb{L}}$, which all give rise to the same basis.)
\item 
For any given $w \in \Waff$ there exists $N(w) \in \Z_{\geq 0}$ such that ${}^q \hspace{-1pt} \uH_w = \uH_w$ for all prime numbers $q \geq N(w)$, but computing explicitly a suitable $N(w)$ is an extremely difficult problem, which is not understood in general in any way at this stage.  Moreover, this integer really depends on $w$: there does not exist a bound that works for all $w$'s, unless $\G$ is a torus.
More generally, although many tools (in particular, computational ones) have been developed to study the $q$-canonical bases, it would be desirable to have more control over its behavior.
 \end{enumerate}
 
As we will explain below, the $p$-canonical basis of $\Haff$ is the appropriate combinatorial tool to describe character formulas for $\G$-modules.

%-----------------------------------------------------------------------
\subsection{Tilting character formula}
\label{ss:tilting-char}
%-----------------------------------------------------------------------

We are finally in a position to state the character formula for indecomposable tilting $\G$-modules.

\begin{thm}
\label{thm:tilting-char-formula}
Let $J \subset \Saff$ be a finitary subset such that $C_J \neq \varnothing$. For any $y,w \in \Waff^J$ we have
\[
d^J_{y,w} = \sum_{z \in W} (-1)^{\ell(z)} \cdot {}^p \hspace{-1pt} h_{zyw_J, ww_J}(1).
\]
\end{thm}

The formula in Theorem~\ref{thm:tilting-char-formula} is of course consistent with~\eqref{eqn:formulas-c-d-translation}. It was first conjectured by Williamson and the second author in~\cite{rw-tilting}. The alternating sum appearing in the right-hand side can be interpreted as an ``antispherical'' $p$-Kazhdan--Lusztig polynomial, i.e.~the coefficient of an element of the $p$-canonical basis of a certain module for $\Haff$ (the antispherical module) rather that of $\Haff$ itself. This formula was inspired by an earlier conjectural formula for these multiplicities when $p \geq h$, due to Andersen~\cite{andersen-tilting} (itself inspired by work of Soergel on the similar question for Lusztig's quantum groups at a root of unity); the latter formula involved ordinary Kazhdan--Lusztig polynomials rather than their $p$-versions, and imposed a bound on $w$.

\begin{remark}
Theorem~\ref{thm:tilting-char-formula} implies a weak form of the ``independence of $p$'' idea explained in
Remark~\ref{rmk:coeff-c-d}\eqref{it:multiplicities-indep-p} for tilting modules.
%Recall the point of view introduced in Remark~\ref{rmk:gps}\eqref{it:gps-fields}.
More specifically, for $y,w \in \Waff^J$, as mentioned in~\S\ref{ss:Haff}, if $p \gg 0$ we have ${}^p \hspace{-1pt} h_{zyw_J, ww_J} = h_{zyw_J, ww_J}$ for all $z \in W$, so that $d^J_{y,w}$ does not depend on $p$.
In fact, in this case the formula in Theorem~\ref{thm:tilting-char-formula} coincides with the one in Andersen's conjecture. However, this does \emph{not} mean that Theorem~\ref{thm:tilting-char-formula} proves Andersen's conjecture in large characteristics, because the bound on $w$ in that conjecture also depends on $p$, so that the number of elements to consider grows when $p$ grows. We currently do not have any understanding of this phenomenon.
\end{remark}

In~\cite{rw-tilting}, the character formula was mainly studied when 
%the special case $J=\varnothing$, assuming that 
$p>h$. (As explained above, in this case one can assume that $J=\varnothing$.) In particular, it was proved that in this setting the formula would follow from another conjecture, of categorical nature, postulating the existence of an action of an appropriate Hecke category $\mathscr{H}_{\mathrm{aff}}^{\Bbbk}$ on $\Rep_{\Waff \bullet 0}(\G)$. This categorical conjecture was also proved in case $\G=\mathrm{GL}_n(\Bbbk)$ for some $n$, proving Theorem~\ref{thm:tilting-char-formula} in this case.

The first proof of Theorem~\ref{thm:tilting-char-formula} for a general reductive group $\G$, but again assuming that $p>h$, was found in joint work with Makisumi and Williamson~\cite{amrw}, building on our earlier work~\cite{ar-reductive}. The strategy used there was inspired by work of Arkhipov--Bezrukavnikov--Ginzburg and Bezrukavnikov--Yun which could be used to provide an alternative proof of Soergel's formula for characters of tilting modules for quantum groups. We will not explain all the details of this approach, but will review some constructions that will be relevant for Section~\ref{sec:support}. Let $G = \G^{(1)}$ be the Frobenius twist of $\G$, and $B = \B^{(1)}$ be its Borel subgroup determined by $\B$. We consider the ``Springer resolution'' $\Spr = G \times^B \mathfrak{n}$, where $\mathfrak{n}$ is the Lie algebra of the unipotent radical of $B$. This variety admits a natural action of $G$; we also consider the action of the multiplicative group $\Gm$ induced by the action on $\mathfrak{n}$ where $z \in \Gm$ acts by multiplication by $z^{-2}$, and the category $\Coh^{G \times \Gm}(\Spr)$ of $(G \times \Gm)$-equivariant coherent sheaves on $\Spr$. In~\cite{ar-reductive}, assuming that $\G$ has simply connected derived subgroup and that $p>h$, we construct a functor
\[
F : \Db \Coh^{G \times \Gm}(\Spr) \to \Db \Rep_{\Wext \bullet 0}(\G)
\]
which is a ``degrading functor'' in the sense that we have a canonical isomorphism of functors $F \circ \langle 1 \rangle [1] \cong F$ (where $\langle 1 \rangle$ is the autoequivalence given by tensoring with the tautological $1$-dimensional $\Gm$-module) which (together with $F$) induces, for any $\mathcal{F},\mathcal{G}$ in $\Db \Coh^{G \times \Gm}(\Spr)$, an isomorphism
\[
\bigoplus_{n \in \Z} \Hom_{\Db \Coh^{G \times \Gm}(\Spr)}(\mathcal{F}, \mathcal{G} \langle n \rangle [n]) \simto \Hom_{\Db \Rep_{\Wext \bullet 0}(\G)}(F(\mathcal{F}), F(\mathcal{G})).
\]
The indecomposable tilting modules in $\Rep_{\Wext \bullet 0}(\G)$ belong to the essential image of this functor, and this allows us to convert the character formula problem into a question in the category $\Db \Coh^{G \times \Gm}(\Spr)$. Using earlier work of Rider and the first author, or equivalently of Mautner and the second author, this question can be further converted into a question regarding mixed perverse sheaves on the affine Grassmannian of the connected reductive group dual to $G$ in the sense of Langlands; the latter question is solved in~\cite{amrw} using a ``Koszul duality'' construction.

Later, two independent proofs of the ``categorical'' conjecture from~\cite{rw-tilting} mentioned above were found by Ciappara~\cite{ciappara} and by Bezrukavnikov and the second author~\cite{br-Hecke}. These results provided new proofs of Theorem~\ref{thm:tilting-char-formula} when $p>h$.
Finally, a proof of Theorem~\ref{thm:tilting-char-formula} in full generality was found by Williamson and the second author in~\cite{rw-smith}. The methods used there are completely different; they rely on an application of ``Smith--Treumann theory'' (a modern reformulation, due to Treumann, of ideas of Smith in algebraic topology) in the context of the geometric Satake equivalence. For a review of this proof, see~\cite{riche-survey}.

%-----------------------------------------------------------------------
\subsection{Characters of simple modules: Lusztig's formula}
\label{ss:characters-simples}
%-----------------------------------------------------------------------

As mentioned in~\S\ref{ss:ind-simp-tilt}, from a character formula for indecomposable tilting modules one can in theory derive a character formula for simple modules, i.e.~compute the coefficients $c^J_{y,w}$. Theorem~\ref{thm:tilting-char-formula} therefore also solves this problem in some sense. However, this answer is very far from explicit, and it would be desirable to have a more concrete formula. 

The approach to this problem was greatly influenced by the following celebrated conjecture, due to Lusztig~\cite{lusztig-pbs}.

\begin{conj}
\label{conj:lusztig}
Assume that $p \geq h$. For any $w \in \Waff^\varnothing$ such that $\langle w \bullet 0 + \rho, \alpha^\vee \rangle \leq p(p-h+2)$ for all $\alpha^\vee \in \Phi^\vee_+$, we have
%in $[\Rep(\G)]$ we have
\[
%[\simp(w \bullet 0)] = \sum_{y \in \Waff^\varnothing} 
c_{y,w}^\varnothing = (-1)^{\ell(w)+\ell(y)} h_{w_\circ y, w_\circ w}(1) \quad \text{for any $y \in \Waff^\varnothing$.}
%\cdot [\coweyl(y \bullet 0)]
\]
\end{conj}

Assuming that $p \geq 2h-3$, all the elements $w \in \Waff^\varnothing$ such that $w \bullet 0$ is restricted satisfy the condition in Conjecture~\ref{conj:lusztig}; in this case, if this conjecture is true, from this formula one can deduce the characters of all simple $\G$-modules using the comments in Remark~\ref{rmk:coeff-c-d}\eqref{it:tensor-product-formulas} and~\eqref{eqn:formulas-c-d-translation}. This conjecture is already discussed at length in several sources, including~\cite{jantzen-char, williamson-alg},~\cite[\S II.8.22]{jantzen} or~\cite[Chap.~I, \S 4]{riche}, so we will not discuss it in detail here. Let us simply note that it is now proved under the assumption that $p$ is larger than a certain bound depending on $\Phi$; a first form of this statement, which did not lead to an explicit formula for the bound, was first proved\footnote{This fact now has several alternative proofs, in particular in~\cite{ar-reductive} or as a consequence of Theorem~\ref{thm:tilting-char-formula}.} by the combination of works of Kazhdan--Lusztig, Tanisaki--Kashiwara and Andersen--Jantzen--Soergel (following a program proposed by Lusztig), and a later refinement of these methods by Fiebig provided an explicit bound on $p$ over which the conjecture holds. (This bound is however difficult to compute in practice, and in any case huge compared with $h$.) Finally, examples constructed by Williamson~\cite{williamson-torsion} showed that no polynomial in $h$ can be a sufficient bound for the validity of the conjecture.

\begin{remark}
Instead of trying to express characters of simple $\G$-modules as linear combinations of characters of induced modules, one can use the fact that the simple $\G$-modules whose highest weight is restricted remain simple as modules over $\G_1\T$ (the preimage of $\T^{(1)} \subset \G^{(1)}$ under the Frobenius morphism of $\G$), and try to express their characters as linear combinations of characters of the baby Verma $\G_1\T$-modules (a family of modules whose characters are easy to describe). An equivalent problem is to describe multiplicities of simple $\G_1\T$-modules in baby Verma modules. It is known that Conjecture~\ref{conj:lusztig} is equivalent to another conjecture of Lusztig, describing these multiplicities in terms of periodic Kazhdan--Lusztig polynomials. In~\cite{rw-simple}, Williamson and the second author define $p$-versions of these periodic polynomials, and deduce from Theorem~\ref{thm:tilting-char-formula} the version of this statement where periodic Kazhdan--Lusztig polynomials are replaced by their $p$-versions, under the assumption that $p>2h-1$. (The bound appearing here is closely related to that appearing for Donkin's formula in Remark~\ref{rmk:coeff-c-d}\eqref{it:tensor-product-formulas}. In fact, using recent advances on this subject it can be improved to $p \geq \max(h,2h-4)$.)
\end{remark}

%-----------------------------------------------------------------------
\subsection{The Finkelberg--Mirkovi{\'c} conjecture}
\label{ss:FM-conj}
%-----------------------------------------------------------------------

In the remainder of this section we outline a slightly different perspective on the question of computing characters of simple modules, offered by a (generalization of a) conjecture due to Finkelberg--Mirkovi{\'c}, which has the merit of providing (less explicit) character formulas which only require that $p>h+1$.

From now on we assume that $\Bbbk$ is an algebraic closure of $\mathbb{F}_p$.
Let $\F$ be an algebraically closed field whose characteristic is positive but different from $p$, and let $H$ be a connected reductive algebraic group over $\F$ which is Langlands dual to $G$. This means in particular that we are given a maximal torus $T_H \subset H$ whose lattice of morphisms to the multiplicative group $\Gm$ (now over $\F$) identifies with $\Y$. Note that the Weyl groups of $\G$ and $H$ identify canonically, so that $\Waff$ is also the semidirect product of the Weyl group of $(H,T_H)$ with the lattice of morphisms from $\Gm$ to $T_H$. We can then consider the loop group $\Loop H$ of $H$, its arc group $\Loop^+ H$, and its affine Grassmannian $\Gr = \Loop^+ H \backslash \Loop H$, with its action of $\Loop H$ induced by multiplication on the right. Below we will in fact only consider the connected component $\Gr^\circ$ of $\Gr$ containing the base point. 

To any finitary subset $J \subset \Saff$ one can associate as in~\cite[\S 3.4]{ar-steinberg} a certain pro-unipotent subgroup scheme $\Iwu^J \subset \Loop H$ (whose action on $\Gr$ stabilizes $\Gr^\circ$) and a morphism of group schemes $\psi_J$ from $\Iwu^J$ to the additive group $\Ga$ over $\F$. (For instance, in case $J=\varnothing$ the group $\Iwu^J$ is the pro-unipotent radical of an Iwahori subgroup and $\psi_J$ is trivial; in case $J=S$ the group $\Iwu^J$ is the pro-unipotent radical of the opposite Iwahori subgroup and $\psi_J$ is induced by a generic character of a maximal connected unipotent subgroup of $H$.) Next, one can choose an Artin--Schreier local system $\mathscr{L}$ on $\Ga$, and consider the derived category $\Db_{\Whit,J}(\Gr^\circ,\Bbbk)$ of bounded constructible complexes of \'etale $\Bbbk$-sheaves on $\Gr^\circ$ which are $(\Iwu^J, \psi_J^*(\mathscr{L}))$-equivariant. (Here, ``$\Whit$'' stands for Whittaker, since this construction is a geometric analogue of Whittaker models in the representation theory of $p$-adic groups.) This category admits a canonical perverse t-structure, whose heart will be denoted $\Perv_{\Whit,J}(\Gr^\circ,\Bbbk)$. The orbits of $\Iwu^J$ on $\Gr^\circ$ are naturally labelled by $W \backslash \Waff$, which is in a natural bijection with $\Waff^\varnothing$, but only those corresponding to elements in $\Waff^J$ support a nonzero $(\Iwu^J, \psi_J^*(\mathscr{L}))$-equivariant local system. Taking $!$-extension, $*$-extension, and intermediate extension of the perversely shifted rank-$1$ equivariant local system on such orbits, one obtains for any $w \in \Waff^J$ a standard perverse sheaf $\Delta^J_w$, a costandard perverse sheaf $\nabla^J_w$, and an intersection cohomology complex $\IC^J_w$.

Given finitary subsets $J \subset K \subset \Saff$, one can define ``averaging functors''
\begin{equation}
\label{eqn:averaging-Whit}
\Perv_{\Whit,J}(\Gr^\circ,\Bbbk) \to \Perv_{\Whit,K}(\Gr^\circ,\Bbbk), \hspace{0.2cm} \Perv_{\Whit,K}(\Gr^\circ,\Bbbk) \to \Perv_{\Whit,J}(\Gr^\circ,\Bbbk).
\end{equation}
(In case $J=\varnothing$, the construction is explained in~\cite[\S 3.6]{ar-steinberg}. The general case is similar.)

The geometric Satake equivalence~\cite{mv} provides an equivalence of monoidal categories between the category $\Perv_{\Loop^+ H}(\Gr^\circ, \Bbbk)$ of $\Loop^+ H$-equivariant perverse sheaves on $\Gr^\circ$ and the category $\Rep(G/\mathrm{Z}(G))$ of finite-dimensional algebraic representations of the quotient of $G$ by its (scheme-theoretic) center. For any finitary subset $J \subset \Saff$, we have an action of $\Perv_{\Loop^+ H}(\Gr^\circ, \Bbbk)$ on $\Perv_{\Whit,J}(\Gr^\circ,\Bbbk)$ by convolution on the left. We also have an action of $\Rep(G/\mathrm{Z}(G))$ on $\Rep(\G)$ where a representation $V$ acts by tensor product with its pullback under the composition $\G \xrightarrow{\Fr} G \to G/\mathrm{Z}(G)$ (where $\Fr$ is the Frobenius morphism), and it follows from the Steinberg tensor product theorem that this action stablizes any ``block'' $\Rep_c(\G)$ with $c \in \X / (\Waff,\bullet)$.

The following conjecture is an extension to general blocks of a conjecture due to Finkelberg--Mirkovi{\'c}~\cite{fm} (which only considered the case $J=\varnothing$ when $C_\varnothing \neq \varnothing$, i.e.~when $p \geq h$).

\begin{conj}
\label{conj:fm}
For any finitary subset $J \subset \Saff$ and any $\lambda \in C_J$ there exists an equivalence of categories
\[
\Psi_\lambda : \Perv_{\Whit,J}(\Gr^\circ,\Bbbk) \simto \Rep_{\Waff \bullet \lambda}(\G)
\]
which satisfies for any $w \in \Waff^J$
\[
\Psi_\lambda(\Delta^J_w) \cong \weyl(w \bullet \lambda), \quad 
\Psi_\lambda(\nabla^J_w) \cong \coweyl(w \bullet \lambda), \quad 
\Psi_\lambda(\IC^J_w) \cong \simp(w \bullet \lambda)
\]
and which intertwines the actions of $\Perv_{\Loop^+ H}(\Gr^\circ, \Bbbk)$ on $\Perv_{\Whit,J}(\Gr^\circ,\Bbbk)$ and of $\Rep(G/\mathrm{Z}(G))$ on $\Rep_{\Waff \bullet \lambda}(\G)$, where these categories are identified via the geometric Satake equivalence.

Moreover, these equivalences can be chosen in such a way that for finitary subsets $J \subset K \subset \Saff$, and for $\lambda \in C_J$ and $\mu \in C_K$, the translation functors $T_\lambda^\mu$ and $T_\mu^\lambda$ (see~\eqref{eqn:translation-functors}) correspond to the natural averaging functors relating $\Perv_{\Whit,J}(\Gr^\circ,\Bbbk)$ and $\Perv_{\Whit,K}(\Gr^\circ,\Bbbk)$, see~\eqref{eqn:averaging-Whit}.
\end{conj}

In the generality considered here, Conjecture~\ref{conj:fm} is completely open. It is however known in the special case $J=\varnothing$ (which, as explained above, is the context of the original conjecture from~\cite{fm}), assuming that $p>h+1$, due to work of Bezrukavnikov and the second author, see~\cite[\S 11]{br-2real}. In~\cite{ar-steinberg,ar-blocks} we also prove various properties of the categories $\Perv_{\Whit,J}(\Gr^\circ,\Bbbk)$ (without any restriction on $p$) which mirror classical results on blocks of $\Rep(\G)$, and thus provide evidence for the conjecture in general.

\begin{remark}
\label{rmk:combinatorics-Perv-Whit}
Let us emphasize a striking aspect of Conjecture~\ref{conj:fm}: the category $\Perv_{\Whit,J}(\Gr^\circ,\Bbbk)$ makes sense for any $J$ and any $\Bbbk$, even when $C_J = \varnothing$. In particular, if true it would provide a ``regular block'' (corresponding to $J=\varnothing$) for any $p$, even when $p<h$ (i.e.~$C_\varnothing = \varnothing$). The geometric analogues of the properties~\eqref{it:translation-coweyl} and~\eqref{it:translation-simples} of~\S\ref{ss:translation} are known in full generality, see~\cite[Lemma~3.3]{ar-steinberg}, so that we have an analogue in this setting of the left-hand formula in~\eqref{eqn:formulas-c-d-translation}. In particular, this reduces the study of the simple objects in all categories $\Perv_{\Whit,J}(\Gr^\circ,\Bbbk)$ to the case $J=\varnothing$, without any restriction on $p$. This also implies that the combinatorics of simple perverse sheaves in such categories does not depend on the choice of Artin--Schreier local system, since this local system does not play any role in case $J=\varnothing$ (because $\psi_\varnothing$ is constant).
\end{remark}

When true for a given $J$, Conjecture~\ref{conj:fm} leads to a formula for the coefficients $(c_{y,w}^J : y,w \in \Waff^J)$.
%expressing simple objects in $\Rep_{\Waff \bullet \lambda}(\G)$ in terms of induced modules. 
In fact, consider for any $y \in \Waff^J$ the map $\chi_y$ from the Grothendieck group $[\Perv_{\Whit,J}(\Gr^\circ,\Bbbk)]$ of $\Perv_{\Whit,J}(\Gr^\circ,\Bbbk)$ to $\Z$ sending the class of an object $\mathscr{F}$ to $\sum_{i \in \Z} (-1)^{i+\ell(y)} \dim \mathsf{H}^i(\mathscr{F}_y)$, where $\mathscr{F}_y$ is the stalk of $\mathscr{F}$ at any point of the orbit corresponding to $y$. It is clear that $( [\Delta_w^J] : w \in \Waff^J )$ is a basis of $[\Perv_{\Whit,J}(\Gr^\circ,\Bbbk)]$, and that $\chi_y(\Delta_w^J)$ equals $1$ if $w=y$, and $0$ otherwise. Hence for any $\mathscr{F}$ we have
\[
[\mathscr{F}] = \sum_{y \in \Waff^J} \chi_y(\mathscr{F}) \cdot [\Delta_y^J].
\]
In particular, if Conjecture~\ref{conj:fm} holds for given $J$, $\lambda$, then for any 
%$\lambda \in C_J$ 
$y,w \in \Waff^J$ we have
%\[
%[\simp(w \bullet \lambda)] = \sum_{y \in \Waff^J} 
\[
c_{y,w}^J = \chi_y(\IC^J_w).
\]
% \cdot [\coweyl(y \bullet \lambda)] \quad \text{
% for any $w \in \Waff^J$.
%\]
Unfortunately the integers $\chi_y(\IC^J_w)$ are not known in general. Essentially the only thing that is known\footnote{To justify this fact, by Remark~\ref{rmk:combinatorics-Perv-Whit} one can assume that $J=\varnothing$, and then use standard arguments to reduce the claim to a similar question for perverse sheaves for the analytic topology on the complex version of $\Gr^\circ$. In this setting it is known that for $p \gg 0$ the graded dimensions of the stalks of $\IC^\varnothing_w$ are the same as for the analogous object with characteristic-$0$ coefficients, which (thanks to results of Kazhdan--Lusztig) are computed by Kazhdan--Lusztig polynomials.} is that for a given $w \in \Waff^J$, if $p$ is large we have
\[
\chi_y(\IC^J_w) = (-1)^{\ell(y)+\ell(w)} \sum_{z \in W_J} (-1)^{\ell(z)} h_{w_\circ yz, w_\circ w}(1)
\]
for any $y \in \Waff^J$,
in accordance with~\eqref{eqn:formulas-c-d-translation} and Conjecture~\ref{conj:lusztig}.

%%%%%%%%%%%%%%%%%%%%%%%%%%%%%%%%
\section{Support varieties}
\label{sec:support}
%%%%%%%%%%%%%%%%%%%%%%%%%%%%%%%%

%-----------------------------------------------------------------------
\subsection{The Frobenius kernel and cohomology}
\label{ss:frob-kernel}
%-----------------------------------------------------------------------

As in~\S\ref{ss:tilting-char} we consider
the Frobenius twist $G = \G^{(1)}$ of $\G$, and let $\Fr: \G \to G$ be the Frobenius morphism, a surjective homomorphism of group schemes.  Let $\G_1 := \ker (\Fr)$ be its scheme-theoretic kernel; this is a finite group scheme over $\Bbbk$, i.e., an affine $\Bbbk$-group scheme corresponding to a finite-dimensional Hopf algebra over $\Bbbk$.  Let $\Rep(\G_1)$ be the category of finite-dimensional algebraic representations of $\G_1$.  For any two $\G_1$-modules $M,N \in \Rep(\G_1)$, one can consider the graded $\Bbbk$-vector space
\[
\Ext^\bullet_{\Rep(\G_1)}(M,N).
\]
If $M$ and $N$ start out as $\G$-modules, then the $\Ext$-groups above retain a residual action of $\G/\G_1 \cong G$.

The $\G_1$-cohomology of a $\G_1$-module $M$, denoted by $\coh^\bullet(\G_1;M)$, is defined by
\[
\coh^\bullet(\G_1;M) := \Ext^\bullet_{\Rep(\G_1)}(\Bbbk,M),
\]
where on the right-hand side $\Bbbk$ is the trivial $\G_1$-module.  We can also regard $\Bbbk$ as the trivial $\G$-module, so if $M$ is a $\G$-module, then $\coh^\bullet(\G_1;M)$ is a graded $G$-module.

In particular, $\coh^\bullet(\G_1;\Bbbk)$ has the structure of a graded $\Bbbk$-algebra with a compatible action of $G$, and it acts naturally on the graded vector space $\Ext^\bullet_{\Rep(\G_1)}(M,N)$ for $M, N$ as above; see~\cite[\S 2.2]{npv}. One is therefore led to ask: what is this ring?  The answer, which we will now recall, was found (assuming $p>h$) in the early 1980s by Andersen--Jantzen~\cite{andersen-jantzen} and Friedlander--Parshall~\cite{friedlander-parshall}.

Let $\fg$ be the Lie algebra of $G$, and let $\cN \subset \fg$ be the variety of nilpotent elements in $\fg$.  Let $\Bbbk[\fg]$ and $\Bbbk[\cN]$ denote the coordinate rings of these varieties.  Of course, $\Bbbk[\fg]$ is the symmetric algebra on $\fg^*$.  
Let the multiplicative group $\Gm$ act on $\fg$ by $z \cdot x = z^{-2}x$ for $z \in \Gm$ and $x \in \fg$.  This $\Gm$-action preserves $\cN$, and it makes the rings $\Bbbk[\fg]$ and $\Bbbk[\cN]$ into graded rings, with the grading determined by the fact that elements of $\fg^*$ (and their images in $\Bbbk[\cN]$) have degree $2$.
The adjoint action of $G$ on $\fg$ commutes with that of $\Gm$, and so the induced actions of $G$ on $\Bbbk[\fg]$ and on $\Bbbk[\cN]$ respect the gradings on those rings.

The results of~\cite{andersen-jantzen, friedlander-parshall} state that when $p > h$, there is a $G$-equivariant graded ring isomorphism
\begin{equation}
\label{eqn:andersen-jantzen}
\coh^\bullet(\G_1;\Bbbk) \cong \Bbbk[\cN].
\end{equation}
(See also~\cite[Chap.~II.12]{jantzen} for a discussion of this result.)
As a consequence:
\begin{itemize}
\item 
for any $M, N \in \Rep(\G_1)$, the graded vector space $\Ext^\bullet_{\Rep(\G_1)}(M,N)$ has the structure of a graded $\Bbbk[\cN]$-module, or equivalently, a $\Gm$-equivariant quasicoherent sheaf on $\cN$;
\item 
if $M$ and $N$ are $\G$-modules, then $\Ext^\bullet_{\Rep(\G_1)}(M,N)$ is a $G$-module, and thus can be regarded as a $(G \times \Gm)$-equivariant quasicoherent sheaf on $\cN$.
\end{itemize}
Thanks to a theorem of Friedlander--Suslin~\cite{friedlander-suslin}, these $\Ext$-groups are finitely generated as modules over $\coh^\bullet(\G_1;\Bbbk)$; i.e., the quasicoherent sheaves mentioned above are in fact coherent.

Describing these coherent sheaves explicitly can be difficult, so one usually first concentrates on the study of their support. Namely,
%The \emph{Humphreys conjecture} is a prediction about the supports of these coherent sheaves in certain cases.  
given a $\G_1$-module $M$, we set
\begin{align*}
V_{\G_1}(M) &:= \supp \bigl( \Ext^\bullet_{\Rep(\G_1)}(M,M) \bigr), \\
\oV_{\G_1}(M) &:= \supp \bigl( \Ext^\bullet_{\Rep(\G_1)}(\Bbbk, M) \bigr) = \supp \bigl( \coh^\bullet(\G_1;M) \bigr).
\end{align*}
These are both (Zariski-)closed $\Gm$-invariant subsets of $\cN$.  We call $V_{\G_1}(M)$ the \emph{support variety} of $M$, and $\oV_{\G_1}(M)$ its \emph{relative support variety}.  If $M$ is a $\G$-module, then $V_{\G_1}(M)$ and $\oV_{\G_1}(M)$ are both $G$-invariant.
% closed subsets of~$\cN$.
For results regarding these supports in case $M$ is an induced or a simple module, see~\cite{npv}. Here we will concentrate on the case $M$ is an indecomposable tilting module, which is the setting for the \emph{Humphreys conjecture}.

\begin{remark}
In case $p \leq h$, the precise structure of $\coh^\bullet(\G_1;\Bbbk)$ is unknown. See~\cite[\S 6]{andersen-jantzen} for some examples. Applying a general result of Suslin--Friedlander--Bendel, one obtains that the spectrum of the even part of this graded ring is homeomorphic to the variety of $p$-nilpotent elements in $\fg$, i.e.~elements annihilated by the restricted $p$th power operation; see~\cite[\S 2.2.10]{npv} for a discussion. Some of the questions considered below make sense in this setting, but they have not been studied as far as we know.
\end{remark}

%-----------------------------------------------------------------------
\subsection{Cells, nilpotent orbits, and Humphreys-style conjectures}
\label{ss:humphreys-conj}
%-----------------------------------------------------------------------

%..........................................................................
\subsubsection{Two-sided cells}
\label{sss:cells}
%..........................................................................

Recall the affine Hecke algebra $\Haff$ defined in~\S\ref{ss:Haff}.  Given a $\Z[v,v^{-1}]$-submodule $A \subset \Haff$, let us say that $A$ is a \emph{based submodule} if the set $A \cap \{\uH_w : w \in \Waff \}$ is a $\Z[v,v^{-1}]$-basis for $A$.  A \emph{based two-sided ideal} is a based submodule that is also a two-sided ideal in $\Haff$.  Define an equivalence relation $\sim_\LR$ on $\Waff$ by
\begin{multline*}
x \sim_\LR y 
\quad \Longleftrightarrow \\
\left(\begin{array}{c}
\text{the smallest based two-sided}\\
\text{ideal containing $\uH_x$}
\end{array}\right) =
\left(\begin{array}{c}
\text{the smallest based two-sided}\\
\text{ideal containing $\uH_y$}
\end{array}\right).
\end{multline*}
This can be extended to an equivalence relation on $\Wext = \Omega \ltimes \Waff$ (see Remark~\ref{rmk:wext-omega}) by declaring that for $x,x' \in \Waff$ and $\omega,\omega' \in \Omega$, we have $\omega x \sim_\LR \omega' x'$ if and only if $x\sim_\LR x'$ in the sense above.  Equivalence classes for this relation are called \emph{two-sided cells}.  An important theorem of Lusztig~\cite{lusztig-cawg4} states that there is a natural bijection
\begin{equation}\label{eqn:lusztig-bij}
\{ \text{two-sided cells in $\Wext$} \} \overset{\sim}{\longleftrightarrow} \{ \text{$G$-orbits in $\cN$} \}.
\end{equation}
With this in mind, we are ready to write down a number of related conjectures, some of which are now theorems.  A brief summary of past results and the current status of these conjectures are given in Tables~\ref{tab:past} and~\ref{tab:current} at the end of the section.

\begin{remark}
\phantomsection
\label{rmk:cells}
\begin{enumerate}
\item
Replacing the Kazhdan--Lusztig basis in the definition above by the $q$-canonical basis of~\S\ref{ss:Haff} (for a prime number $q$) one obtains the notion of two-sided $q$-cells in $\Waff$. This notion was investigated, for a general crystallographic Coxeter system, by L.~T.~Jensen in his PhD thesis. See~\cite{ahr-conj} for some conjectures on these cells in the specific setting of affine Weyl groups, which are closely related to the topic of this section. 
\item
The connection between support varieties and two-sided cells becomes clearer after the introduction of two other notions that are used to fill the gap between them, namely antispherical right ($q$-)cells and weight cells. See~\cite{ahr-conj} for a discussion of these topics and the relation with the questions discussed here. (Weight cells are closely related with the tensor ideals discussed in~\S\ref{ss:tensor-ideals} below.)
\item
Lusztig's bijection in~\cite{lusztig-cawg4} is stated in terms of nilpotent orbits for the \emph{complex} reductive algebraic group corresponding to $G$. However nilpotent orbits for this group are in a canonical bijection with those of $G$ (under mild assumptions on $p$), see~\S\ref{sss:vecbdle} below.
\end{enumerate}
\end{remark}

%..........................................................................
\subsubsection{Traditional and relative Humphreys conjectures for reductive groups}
\label{sss:humphreys}
%..........................................................................

From now on we assume that $p>h$. We will denote by $\Wext^\varnothing \subset \Wext$ the subset consisting of element $w$ which are of minimal length in their coset $Ww$; in fact we have $\Wext^\varnothing = \bigsqcup_{\omega \in \Omega} \Waff^\varnothing \cdot \omega$ using the notation of~\S\ref{ss:translation}.
The following statement was originally proposed by Humphreys in~\cite[\S 12]{humphreys}.

\begin{conj}[Traditional Humphreys conjecture]
\label{conj:humphreys}
Let $w \in \Wext^\varnothing$, and let $C$ be the $G$-orbit in $\cN$ corresponding under~\eqref{eqn:lusztig-bij} to the two-sided cell containing $w$.  Then
\[
V_{\G_1}(\til(w \bullet 0)) = \overline{C}.
\]
\end{conj}

\begin{remark}
\begin{enumerate}
\item 
Although the preceding statement is widely known as ``Hum\-phreys' conjecture,'' Hum\-phreys did not consider himself to have ``conjectured'' it to be true, but rather to have merely proposed it as a statement worthy of further inquiry.  In~\cite{humphreys}, the statement is labelled as a ``Hypothesis.''
\item The ``hypothesis'' in~\cite{humphreys} describes $V_{\G_1}(\til(\lambda))$ for all $\lambda \in \X_+$ which are ``regular,'' i.e.~have trivial stabilizer in $\Waff$, and not just for $\lambda \in \Wext^\varnothing \bullet 0$.  But as explained in~\cite[\S8.2]{ahr}, the statement in Conjecture~\ref{conj:humphreys} actually determines $V_{\G_1}(\til(\lambda))$ for all $\lambda \in \X_+$, so it is equivalent to Humphreys' original version.
\end{enumerate}
\end{remark}

We will now state several variants of Conjecture~\ref{conj:humphreys} which were not explicitly considered by Humphreys, but to which his name is usually attached too. The first one is concerned with relative support varieties. For this we will denote by $\Wext^{\varnothing,\prime} \subset \Wext$ the subset of elements $w$ which are of minimal length in their double coset $WwW$. (Note that $\Wext^{\varnothing,\prime} \cap \Waff$ does \emph{not} coincide with the subset $\Waff^S$ introduced in~\S\ref{ss:translation}; rather, it strictly contains it unless $\G$ is a torus.) It is easy to show that
\[
\coh^\bullet(\G_1; \til(\lambda)) = 0
\qquad\text{if}\qquad
\lambda \notin \Wext^{\varnothing,\prime} \bullet 0.
\]
(See, for instance,~\cite[Lemma~8.7]{ahr}.)  In other words, our conjecture needs only deal with $\til(w \bullet 0)$ for $w \in \Wext^{\varnothing,\prime}$.  Next, we have natural bijections
\[
\X_+ \overset{\sim}{\to} W \backslash \Wext / W \overset{\sim}{\leftarrow} \Wext^{\varnothing,\prime}.
\]
For $\lambda \in \X_+$, let $w_\lambda$ be the corresponding element of $\Wext^{\varnothing,\prime}$, or in other words the unique element of minimal length in $Wt_\lambda W$.  (It will be convenient later in this section to label elements of $\Wext^{\varnothing,\prime}$ by dominant weights in this way.)

\begin{conj}[Relative Humphreys conjecture]\label{conj:rel-humphreys}
Let $\lambda \in \X_+$, and let $C$ be the $G$-orbit in $\cN$ corresponding under~\eqref{eqn:lusztig-bij} to the two-sided cell containing $w_\lambda$.  Then
\[
\oV_{\G_1}(\til(w_\lambda \bullet 0)) = \overline{C}.
\]
\end{conj}

%..........................................................................
\subsubsection{Scheme-theoretic versions}
\label{sss:Humphreys-scheme-theoretic}
%..........................................................................

%We now return to the standard 
In~\S\ref{ss:frob-kernel}, $V_{\G_1}(M)$ and $\oV_{\G_1}(M)$ are defined as closed \emph{sets} in $\cN$, but one can also take into account \emph{schematic} structures.  The following statements are similar in spirit to, but a priori stronger than, those in~\S\ref{sss:humphreys}.

\begin{conj}[Scheme-theoretic traditional Humphreys conjecture]\label{conj:sch-humphreys}
Let $w \in \Wext^\varnothing$, and let $C$ be the $G$-orbit in $\cN$ corresponding under~\eqref{eqn:lusztig-bij} to the two-sided cell containing $w$.  Then the coherent sheaf $\Ext^\bullet_{\Rep(\G_1)}(\til(w \bullet 0),\til(w \bullet 0))$ is scheme-theoretically supported on the (reduced) closed subscheme $\overline{C}$.
\end{conj}

\begin{conj}[Scheme-theoretic relative Humphreys conjecture]\label{conj:sch-rel-humphreys}
Let $\lambda \in \X_+$, and let $C$ be the $G$-orbit in $\cN$ corresponding under~\eqref{eqn:lusztig-bij} to the two-sided cell containing $w_\lambda$.  Then the coherent sheaf $\coh^\bullet(\G_1; \til(w_\lambda \bullet 0))$ is scheme-theoretically supported on the (reduced) closed subscheme $\overline{C}$.
\end{conj}

%..........................................................................
\subsubsection{Quantum group versions}
\label{sss:quantum-versions}
%..........................................................................

In this subsection, 
%we temporarily allow $p$ to be any positive integer, not necessarily prime, and let $\zeta \in \C$ be a primitive $p$th root of unity. 
we assume that $\G$ is semisimple and simply connected. Let $\G_\C$ be the complex semisimple algebraic group with the same root system as $\G$, and let $\fg_\C$ be its Lie algebra. (In this setting, there is no point distinguishing the group from its Frobenius twist, so we also set $G_\C=\G_\C$.) Let $\U_v(\fg_\C)$ be the Lusztig integral form of the quantum group associated with the root system of $\G$ (or $\G_\C$).  This is an algebra over $\Z[v,v^{-1}]$.  Let $l$ be an odd integer, coprime to $3$ if $\Phi$ has a component of type $G_2$, and let $\zeta \in \C$ be a primitive $l$th root of unity. Let
\[
\U_\zeta(\fg_\C) = \C \otimes_{\Z[v,v^{-1}]} \U_v(\fg_\C)
\]
be the specialization obtained by setting $v = \zeta$.  This is a Hopf algebra over $\C$, and its representation theory has a great deal in common with that of $\G$; see e.g.~\cite[Chap.~II.H]{jantzen} for a brief review. For instance, there are quantum induced modules $\coweyl_\zeta(\lambda)$, quantum Weyl modules $\weyl_\zeta(\lambda)$, and quantum indecomposable tilting modules $\til_\zeta(\lambda)$.  There is also a (surjective) ``quantum Frobenius map''
$
\Fr_\zeta: \U_\zeta(\fg_\C) \to U(\fg_\C)$,
where $U(\fg_\C)$ is the universal enveloping algebra of $\fg_\C$, and a finite-dimensional sub-Hopf algebra
\[
\su_\zeta(\fg_\C) \subset \U_\zeta(\fg_\C),
\]
called the \emph{small quantum group}, that can informally be thought of as the ``Hopf-algebra-theoretic kernel'' of $\Fr_\zeta$.  In analogy with~\eqref{eqn:andersen-jantzen}, by results of Ginzburg--Kumar~\cite{ginzburg-kumar}, if $l>h$ we have
\[
\coh^\bullet(\su_\zeta(\fg_\C); \C) \cong \C[\cN_\C]
\]
where $\cN_\C \subset \fg_\C$ is the nilpotent cone.
As a consequence, if $M$ and $N$ are finite-dimensional $\U_\zeta(\fg_\C)$-modules, then $\Ext^\bullet_{\su_\zeta(\fg_\C)}(M,N)$ can be regarded as a $(\G_\C \times \C^\times)$-equivariant quasicoherent sheaf on $\cN_\C$. 

For $w \in \Wext$, one can define the ordinary and relative support varieties
\[
V_{\su_\zeta(\fg_\C)}(\til_\zeta(w \bullet 0))
\qquad\text{and}\qquad
\oV_{\su_\zeta(\fg_\C)}(\til_\zeta(w \bullet 0))
\]
by mimicking the definitions in~\S\ref{ss:frob-kernel}.  
With these definitions in place, one can consider the quantum variants of each of Conjectures~\ref{conj:humphreys}, \ref{conj:rel-humphreys}, \ref{conj:sch-humphreys}, and~\ref{conj:sch-rel-humphreys}.

\begin{table}
\begin{center}
\begin{tabular}{lll}
\textit{Year} & \textit{Authors} & \textit{Theorems} \\
\hline
1998 & Humphreys~\cite{humphreys} & proposed Conjecture~\ref{conj:humphreys} \\
1998 & Ostrik~\cite{ostrik} & quantum traditional for type $A$ \\
2006 & Bezrukavnikov~\cite{bezrukavnikov} & quantum relative\\
&& scheme-theoretic quantum relative\\
2018 & Hardesty~\cite{hardesty} & traditional for type $A$ \\
2019 & Achar--Hardesty--Riche~\cite{ahr} & traditional \& relative for $p \gg 0$ \\
&& quantum traditional \\
2024 & Achar--Hardesty~\cite{ah1} & scheme-theoretic relative for type $A$ \\
2024 & Achar--Hardesty--Riche~\cite{ahr2,ahr3} & relative \\
& $+$ Achar--Hardesty~\cite{ah2} & \\
\end{tabular}
\end{center}
\caption{History of results on the Humphreys conjecture}\label{tab:past}
\end{table}

\begin{table}
\begin{center}
\begin{tabular}{l|c|c}
& \textit{Classical} & \textit{Quantum} \\
\hline
\textit{Traditional} & 
known for type $A$ or for $p \gg 0$ & known \\
\hline
\textit{Relative} & known & known \\
\hline
\begin{tabular}{@{}l@{}}
\textit{Scheme-theoretic} \\ \textit{traditional}
\end{tabular}
& known for type $A$ or for $p \gg 0$ & known \\
\hline
\begin{tabular}{@{}l@{}}
\textit{Scheme-theoretic} \\ \textit{relative}
\end{tabular}
& known for type $A$ or for $p \gg 0$ & known
\end{tabular}
\end{center}
\caption{Summary of current status}\label{tab:current}
\end{table}

%-----------------------------------------------------------------------
\subsection{Cohomology and coherent sheaves on the Springer resolution}
\label{ss:cohomology-coherent}
%-----------------------------------------------------------------------

In the remainder of this section we will discuss the proofs of (some of) the results in Table~\ref{tab:past}.

We now come back to the setting of~\S\S\ref{sss:humphreys}--\ref{sss:Humphreys-scheme-theoretic}, assuming in particular that $p>h$.
Recall from~\S\ref{ss:tilting-char} that we have a degrading functor
\[
F : \Db \Coh^{G \times \Gm}(\Spr) \to \Db \Rep_{\Wext \bullet 0}(\G),
\]
and that the tilting modules in $\Rep_{\Wext \bullet 0}(\G)$ lie in the essential image of this functor.  In fact, by imposing a kind of ``self-duality'' condition, one can pick out a canonical preimage for each tilting module $\til(w \bullet 0)$ with $w \in \Wext^\varnothing$, which we denote by
\[
\til^{\mathrm{gr}}(w \bullet 0) \in \Db \Coh^{G \times \Gm}(\Spr).
\]
(These objects admit an intrinsic characterization, independent of the functor $F$, in the spirit of the definition of the Juteau--Mautner--Williamson parity complexes; see~\cite[proof of Proposition~9.1]{ahr} for details.)
Let $\pi: \Spr \to \cN$ be the map that sends $(g,x) \in G \times^B \mathfrak{n}$ to $\mathrm{Ad}(g)(x) \in \cN$.  This map is proper (and indeed a resolution of singularities of $\cN$), so it gives rise to a (derived) push-forward functor
\[
\pi_*:  \Db \Coh^{G \times \Gm}(\Spr) \to \Db\Coh^{G \times \Gm}(\cN).
\]
For $\lambda \in \X_+$, let
\[
\fS_\lambda := \pi_*\til^{\mathrm{gr}}(w_\lambda \bullet 0)\quad \in \Db\Coh^{G \times \Gm}(\cN).
\]
Then, according to~\cite[Proposition~9.1]{ahr}, there is a $(G \times \Gm)$-equivariant isomorphism of graded $\Bbbk[\cN]$-modules
\[
\coh^\bullet(\G_1; \til(w_\lambda \bullet 0)) = 
%\bigoplus_{i \in \Z} 
R^\bullet\Gamma(\cN, \fS_\lambda).
%\langle -i \rangle.
\]
(For $w \in \Wext^\varnothing \smallsetminus \Wext^{\varnothing,\prime}$, we have $\pi_*\til^{\mathrm{gr}}(w \bullet 0) = 0$, paralleling the fact that $\coh^\bullet(\G_1; \til(w \bullet 0)) = 0$, see~\S\ref{sss:humphreys}.)  Thus, the relative or scheme-theoretic relative Humphreys conjectures can be rephrased as conjectures about the set-theoretic or scheme-theoretic support of $\fS_\lambda$, which constitutes the main tool for its study in the setting of a general reductive group.

%-----------------------------------------------------------------------
\subsection{The Lusztig--Vogan bijection}
%-----------------------------------------------------------------------

We now take a slight detour to discuss some results related to nilpotent orbits.

%..........................................................................
\subsubsection{Parametrizing vector bundles}
\label{sss:vecbdle}
%..........................................................................

Let $C \subset \cN$ be a $G$-orbit.  Choose a (closed) point $x \in C$, and let $G^x$ be its scheme-theoretic stabilizer in $G$. In fact, $G^x$ is in reduced, i.e., this group scheme is actually an algebraic group over $\Bbbk$, see e.g.~the discussion in~\cite[\S 3.1]{ahr3}. There is an equivalence of abelian categories
\begin{equation}
\label{eqn:coh-nilp}
\Coh^G(C) \cong \Rep(G^x)
\end{equation}
between the categories of $G$-equivariant coherent sheaves on $C$ and of finite-dimen\-sional representations of $G^x$.
In particular, every object in $\Coh^G(C)$ admits a composition series, i.e., a finite filtration whose subquotients are simple objects, and the simple objects in $\Coh^G(C)$ correspond under~\eqref{eqn:coh-nilp} to the irreducible $G^x$-modules.  Let
\[
\Sigma_C := \{ \text{isomorphism classes of simple objects in $\Coh^G(C)$} \},
\]
and given an isomorphism class $\sigma \in \Sigma_C$, let $\cL(\sigma)$ be a simple object in this class.

We can transfer other representation-theoretic notions across~\eqref{eqn:coh-nilp}.   For instance, let $(G^x)^\circ$ be the identity connected component of $G^x$,
%.  Then $\Rep(G^x/(G^x)^\circ)$ is identified with a full subcategory of $\Rep(G^x)$.  The corresponding full subcategory on the left-hand side of~\eqref{eqn:coh-nilp} is the category of $G$-equivariant flat vector bundles on $C$.
let $G^x_\unip \subset (G^x)^\circ$ be the unipotent radical of $(G^x)^\circ$, and set
\[
G^x_\red := G^x/ G^x_\unip.
\]
This is a (possibly disconnected) reductive group.  The category $\Rep(G^x_\red)$ identifies with a full subcategory of $\Rep(G^x)$, and this subcategory contains every simple $G^x$-module.

Our assumption that $p > h$ implies that $p \nmid |G^x_\red/(G^x_\red)^\circ| = |G^x/(G^x)^\circ|$; according to~\cite{ahr2}, the categorical structure described in~\S\ref{ss:ind-simp-tilt} therefore generalizes to $\Rep(G^x_\red)$.  In particular, there is a well-defined notion of \emph{tilting module} for $\Rep(G^x_\red)$, and isomorphism classes of indecomposable tilting modules are in a canonical bijection with those of simple modules.
An object of $\Coh^G(C)$ is called a \emph{tilting vector bundle} if it corresponds to a tilting $G^x_\red$-module via~\eqref{eqn:coh-nilp}.  By the comments above there is a canonical bijection between $\Sigma_C$ and the set of isomorphism classes of indecomposable tilting vector bundles; for $\sigma \in \Sigma_C$, we let $\cT(\sigma)$ be an indecomposable tilting vector bundle in the corresponding class.

Finally, let
\[
\Xi_G := \{ (C,\sigma) \mid \text{$C \subset \cN$ a $G$-orbit, and $\sigma \in \Sigma_C$} \}.
\]
It turns out that the set $\Xi_G$ essentially depends only on the root datum of $G$, and not on the field $\Bbbk$.  In more detail, under mild assumptions on $p$, we have:
\begin{enumerate}
\item by the Bala--Carter theorem, the set of $G$-orbits on $\cN$ admits a parametrization that is independent of $\Bbbk$.
\item as shown in~\cite{ahr3}, for $x \in \cN$, the root datum of $(G^x_\red)^\circ$ and the action of the group $G^x_\red/(G^x_\red)^\circ$ on the weight lattice of $(G^x_\red)^\circ$ depend only on the Bala--Carter label of the $G$-orbit containing $x$, and not on $\Bbbk$.
\end{enumerate}
As a consequence, $\Xi_G$ depends only on the root datum of $G$.

%..........................................................................
\subsubsection{The complex case}
\label{sss:vecbdle-C}
%..........................................................................

The constructions from~\S\ref{sss:vecbdle} make sense also in the context of~\S\ref{sss:quantum-versions}.
%for the complex algebraic group $G_\C$.
% with the same root datum as $G$. 
 The counterpart of $\Xi_G$ in this setting will be denoted $\Xi_{G_\C}$; it is also in a canonical bijection with $\Xi_G$. In this setting, around 1990, Lusztig~\cite{lusztig-cawg4} and Vogan~\cite{vogan} independently predicted the following statement.

\begin{conj}\label{conj:lv}
There is a ``natural'' bijection $\X_+ \overset{\sim}{\longleftrightarrow} \Xi_{G_\C}$.
\end{conj}

To make this meaningful, we must clarify the meaning of ``natural'': at a minimum, we must specify some properties that the bijection is expected to have.  In Vogan's version, the expected properties are formulated in terms of associated varieties of Harish-Chandra modules for complex Lie groups.  In Lusztig's version (which is a priori different from Vogan's), the bijection is required to be compatible with~\eqref{eqn:lusztig-bij}, in the following sense:
\begin{multline}\label{eqn:lv-compat}
\left(
\begin{array}{c}
\text{$\lambda \in \X_+$ corresponds to $(C,\sigma) \in \Xi_{G_\C}$}\\
\text{under the bijection of Conjecture~\ref{conj:lv}}
\end{array}
\right)
\\
\text{implies that}\quad
\left(
\begin{array}{c}
\text{the two-sided cell containing $w_\lambda$}\\
\text{corresponds to $C$ under~\eqref{eqn:lusztig-bij}}
\end{array}
\right).
\end{multline}

Vogan's version of Conjecture~\ref{conj:lv} was proved by the first author for $G_\C = \mathrm{GL}_n(\C)$ in~\cite{a-thesis, a-equivkthy}, and then by Bezrukavnikov~\cite{bez-qes} for all $G_\C$.  Subsequently, Bezrukavnikov showed~\cite{bez-tcac, bez-psaf} that this bijection also satisfies the desiderata in Lusztig's version of the conjecture, so that there is indeed a well-defined \emph{Lusztig--Vogan bijection} for complex reductive groups.

%..........................................................................
\subsubsection{The general case}
%..........................................................................

In~\cite{achar-pcoh}, the first author constructed a bijection $\X_+ \overset{\sim}{\longleftrightarrow} \Xi_G$ by adapting Bezrukavnikov's argument from~\cite{bez-qes}.  This bijection could be called the ``Lusztig--Vogan bijection over $\Bbbk$.''  Then, in~\cite{ahr3} it was shown that the triangle
\[
\begin{tikzcd}[row sep=0pt]
& \Xi_G \ar[dd, "\wr"] \\
\X_+ \ar[ur, "\sim"] \ar[dr, "\sim"'] \\
& \Xi_{G_\C}
\end{tikzcd}
\]
commutes.  Here, the diagonal arrows are the $\Bbbk$- and $\C$-versions of the Lusztig--Vogan bijection, and the vertical arrow is the identification discussed in~\S\S\ref{sss:vecbdle}--\ref{sss:vecbdle-C}.  Thus, there is a canonical Lusztig--Vogan bijection for $G$ that depends only on its root datum.

%-----------------------------------------------------------------------
\subsection{Remarks on the proof of the Humphreys conjecture}
%-----------------------------------------------------------------------

The various incarnations of the Humphreys conjecture predict the \emph{supports} of certain coherent sheaves on $\cN$, but one could also ask for more information about what the sheaves themselves are.  In particular, in the relative case, one could ask for more information about the objects $\fS_\lambda$. In the quantum case of~\S\ref{sss:quantum-versions} (where a counterpart of the degrading functor $F$ from~\S\ref{ss:tilting-char} was constructed by Arkhipov--Bezrukavnikov--Ginzburg long before its modular version), these objects were identified by Bezrukavnikov in~\cite{bezrukavnikov}.

\begin{thm}[Bezrukavnikov~\cite{bezrukavnikov}]
\label{thm:bezru}
Suppose $\lambda \in \X_+$ corresponds under the Lusztig--Vogan bijection to $(C,\sigma)$.  Then the complex $\fS_{\lambda,\C}$ is the simple perverse-coherent sheaf $\mathcal{IC}(C,\cL_\C(\sigma))$.
\end{thm}

See~\cite{bez-qes} or~\cite{achar-on} for background and references on the theory of ``perverse-coherent sheaves'' on $\cN$.  By construction, the simple object $\mathcal{IC}(C,\cL_\C(\sigma))$ is supported scheme-theoretically on $\overline{C}$, so Theorem~\ref{thm:bezru} immediately implies the quantum versions of Conjectures~\ref{conj:rel-humphreys} and~\ref{conj:sch-rel-humphreys}.

In~\cite{ahr}, Hardesty and the authors carried out a strategy of ``lifting to characteristic $0$'' to compare the support of $\fS_\lambda$ to that of its complex counterpart.  This lifting can be carried out when the characteristic $p$ of $\Bbbk$ is sufficiently large. In this way we obtained the proof of Conjecture~\ref{conj:rel-humphreys} (and, although it was not explicitly discussed therein, Conjecture~\ref{conj:sch-rel-humphreys} as well) for $p \gg 0$, i.e., $p$ larger than a nonexplicit bound (from the point of view introduced in Remark~\ref{rmk:gps}\eqref{it:gps-fields}).

The same paper also clarifies the relationship between the traditional and relative versions of the Humphreys conjecture: the traditional version implies the relative version~\cite[Remark~9.41(1)]{ahr}; conversely, the relative version together with a technical condition on cells in $\Wext$ implies the traditional version~\cite[Lemma~8.11]{ahr}.  We thus deduce Conjecture~\ref{conj:humphreys} (and, implicitly, Conjecture~\ref{conj:sch-humphreys}) for $p \gg 0$.  Similarly, the quantum versions of these conjectures follow from Theorem~\ref{thm:bezru}.
%~\cite{bezrukavnikov}.

Conceptually, the key to the arguments outlined above is that Theorem~\ref{thm:bezru} identifies the $\fS_\lambda$'s with objects that are constructed in a totally different way, not involving the functors $F$ or $\pi_*$.  To make further progress on the Humphreys conjectures, one might seek such a description of the $\fS_\lambda$'s that is valid over $\Bbbk$.
(Perverse-coherent sheaves make sense over $\Bbbk$, but in general the $\fS_\lambda$'s do \emph{not} coincide with simple perverse-coherent sheaves. Note also that the characterization of the objects $\til^{\mathrm{gr}}(w \bullet 0)$ alluded to in~\S\ref{ss:cohomology-coherent} 
%of course provides a characterization of the objects $\fS_\lambda$, which is however 
is not sufficient for the present discussion.)
  As a step in this direction, in~\cite{ah1}, Hardesty and the first author showed that the $\fS_\lambda$'s are (up to a grading shift) precisely the indecomposable objects in the coheart of a certain co-t-structure on $\Db\Coh^{G \times \Gm}(\cN)$. In the case when $\G=\mathrm{GL}_n(\Bbbk)$, they were also able (using ``parabolic'' versions of the co-t-structure above, and the fact that in this case each nilpotent orbit is a Richardson orbit) to prove Conjectures~\ref{conj:sch-humphreys} and~\ref{conj:sch-rel-humphreys} under the assumption $p>n$. (In this case, Conjecture~\ref{conj:humphreys} had been proved earlier under the same assumption by Hardesty~\cite{hardesty}.)

In~\cite{ah2}, a different construction of the same co-t-structure is given.  This construction is closely tied with the geometry of nilpotent orbits: to each pair $(C,\sigma) \in \Xi_G$ one associates an object
\[
\cS(C,\sigma) \in \Db\Coh^{G \times \Gm}(\cN)
\]
such that the following hold:
\begin{enumerate}
\item the object $\cS(C,\sigma)$ is supported (set-theoretically) on $\overline{C}$;
\item the restriction $\cS(C,\sigma)_{|C}$ is the tilting vector bundle $\cT(\sigma)$;
\item the collection $\{ \cS(C,\sigma) : (C,\sigma) \in \Xi_G \}$ is exactly the set of indecomposable objects (up to grading shift) in the coheart of the co-t-structure.
\end{enumerate}
Thus, each $\cS(C,\sigma)$ must be isomorphic to some $\fS_\lambda$. More specifically, we have the following statement.

\begin{thm}[\cite{ah2}]
Suppose $\lambda \in \X_+$ corresponds under the Lusztig--Vogan bijection to $(C,\sigma)$.  Then $\fS_{\lambda}$ is isomorphic to $\cS(C,\sigma)$.
\end{thm}

When combined with the results from~\cite{ahr,ahr2} discussed above, this yields Conjecture~\ref{conj:rel-humphreys} in full generality.  

\begin{remark}
\phantomsection
\begin{enumerate}
\item The construction of the $\cS(C,\sigma)$'s in~\cite{ah2} does \emph{not} indicate whether their support is reduced, so when $p$ is not large enough for the results of~\cite{ahr} to apply, Conjecture~\ref{conj:sch-rel-humphreys} remains open.
\item Although~\cite[Lemma~8.11]{ahr} gives a method for deducing Conjecture~\ref{conj:humphreys} from Conjecture~\ref{conj:rel-humphreys}, it requires checking a technical condition on cells in $\Wext$ that we do not know how to do.  Thus, Conjecture~\ref{conj:humphreys} and, a fortiori, Conjecture~\ref{conj:sch-humphreys} remain open outside the range covered by~\cite{ahr} or~\cite{hardesty}.
\end{enumerate}
\end{remark}

%-----------------------------------------------------------------------
\subsection{Tensor ideals}
\label{ss:tensor-ideals}
%-----------------------------------------------------------------------

We conclude with some brief comments about a (conjectural) outgrowth of the ideas above.  Recall from~\S\ref{ss:ind-simp-tilt} that the category $\Tilt(\G)$ is closed under tensor product.  A subcategory $\mathcal{I} \subset \Tilt(G)$ is called a \emph{tensor ideal}\footnote{Often, this notion is called a \emph{thick} tensor ideal. Since all the tensor ideals we want to consider are thick, we omit this adjective.} if it is closed under direct sums and direct summands, and if for any $M \in \mathcal{I}$ and any $N \in \Tilt(\G)$ we have $M \otimes N \in \mathcal{I}$.  The \emph{principal tensor ideal} generated by an indecomposable tilting module $M$ is the smallest tensor ideal containing $M$.  

In~\cite{ostrik-tensor}, Ostrik classified the principal tensor ideals in the quantum-group version $\Tilt(\U_\zeta(\fg_\C))$.  When combined with~\eqref{eqn:lusztig-bij}, Ostrik's result can be interpreted as a bijection
\[
\{\text{principal tensor ideals in $\Tilt(\U_\zeta(\fg_\C))$}\} \overset{\sim}{\longleftrightarrow} \{ \text{$G_\C$-orbits in $\cN_\C$} \}.
\]

However, for reductive groups over $\Bbbk$, the category $\Tilt(\G)$ always has infinitely many principal tensor ideals (unless $\G$ is a torus), so any classification theorem for them must be more complicated than Ostrik's statement.  In~\cite{ahr-conj}, the authors (together with Hardesty) proposed a conjectural classification that was motivated by an examination of the data appearing in the Lusztig--Vogan bijection.  Specifically, for each orbit $C \subset \cN$, the category $\Tilt(\Coh^G(C))$ of tilting vector bundles is a tensor category, and one can consider the set
\[
\Upsilon_1 :=  \left\{ (C, \mathcal{I}) \,\Big|\,
\begin{array}{c}
\text{$C \subset \cN$ a $G$-orbit, and} \\
\text{$\mathcal{I} \subset \Tilt(\Coh^G(C))$ a principal tensor ideal}
\end{array}
\right\}.
\]
We conjecture that there is a bijection
\begin{equation}
\label{eqn:ahr-conj}
\{\text{principal tensor ideals in $\Tilt(\G)$}\} \overset{\sim}{\longleftrightarrow} \Upsilon_1.
\end{equation}
Suppose~\eqref{eqn:ahr-conj} is true.  Via~\eqref{eqn:coh-nilp}, we identify $\Tilt(\Coh^G(C))$ with $\Tilt(G^x_\red)$, and then we can apply the bijection~\eqref{eqn:ahr-conj} to the various reductive groups\footnote{A subtlety here is that $G^x_\red$ may be disconnected, so one should first formulate a refinement of~\eqref{eqn:ahr-conj} that accomodates disconnected groups.} $G^x_\red$ to obtain that $\Upsilon_1$ is in bijection with the set
\[
{\small
\Upsilon_2 := \left\{ (C, C', \mathcal{I}) \,\Big|\,
\begin{array}{c}
\text{$C=G \cdot x \subset \cN$ a $G$-orbit, $C' \subset \cN_{(G^x_\red)^{(1)}}$ a $(G^x_\red)^{(1)}$-orbit,} \\
\text{and $\mathcal{I} \subset \Tilt(\Coh^{(G^x_\red)^{(1)}})$ a principal tensor ideal}
\end{array}
\right\}.
}
\]
This process can be repeated to obtain $\Upsilon_3, \Upsilon_4, \ldots$, providing a kind a ``recursive'' parametrization of principal tensor ideals. (See~\cite[\S 6]{ahr-conj} for an explicit description of this process for the group $\mathrm{GL}_n(\Bbbk)$.)

The truth or falsity of this conjecture is closely related to the question of how antispherical right $p$-cells are related to ordinary antispherical right cells (see Remark~\ref{rmk:cells}), and perhaps also (as suggested by Bezrukavnikov) to questions about cohomology and support varieties for higher Frobenius kernels $\G_r := \ker(\Fr^r: \G \to \G^{(r)})$.

%-----------------------------------------------------------------------
\section*{Acknowledgements}
%-----------------------------------------------------------------------

This material is based upon work supported by the National
Science Foundation under Grant Nos.~DMS-1500890, DMS-1802241, and DMS-2202012.
This project has received
funding from the European Research Council (ERC) under the European Union's Horizon 2020
research and innovation programme (grant agreement No.~101002592).

We thank our collaborators involved in the development of the program outlined above. In particular, the contribution of Geordie Williamson was essential for the work presented in Section~\ref{sec:characters}, and all the results presented in Section~\ref{sec:support} were obtained as collaborations with William Hardesty. The influence of Roman Bezrukavnikov on our work, via collaborations, discussions, or study of his papers, also cannot be overestimated.

% SIAM recommends using BibTeX
% if using BibTeX
%\bibliographystyle{siamplain}
%\bibliography{ref_icm}

\begin{thebibliography}{10}

\bibitem{a-thesis}
{\sc P.~N. Achar}, {\em Equivariant coherent sheaves on the nilpotent cone for
  complex reductive {L}ie groups}, {P}h{D} thesis, Massachusetts Institute of
  Technology, 2001.

\bibitem{a-equivkthy}
{\sc P.~N. Achar}, {\em On the equivariant {$K$}-theory of the nilpotent cone
  in the general linear group}, Represent. Theory, 8 (2004), pp.~180--211,
  \url{https://doi.org/10.1090/S1088-4165-04-00243-2}.

\bibitem{achar-pcoh}
{\sc P.~N. Achar}, {\em Perverse coherent sheaves on the nilpotent cone in good
  characteristic}, in Recent developments in {L}ie algebras, groups and
  representation theory, vol.~86 of Proc. Sympos. Pure Math., Amer. Math. Soc.,
  Providence, RI, 2012, pp.~1--23,
  \url{https://doi.org/10.1090/pspum/086/1409}.

\bibitem{achar-on}
{\sc P.~N. Achar}, {\em On exotic and perverse-coherent sheaves}, in
  Representations of reductive groups, vol.~312 of Progr. Math.,
  Birkh\"{a}user/Springer, Cham, 2015, pp.~11--49,
  \url{https://doi.org/10.1007/978-3-319-23443-4\_2}.

\bibitem{ah1}
{\sc P.~N. Achar and W.~Hardesty}, {\em Co-{$t$}-structures on derived
  categories of coherent sheaves and the cohomology of tilting modules},
  Represent. Theory, 28 (2024), pp.~49--89,
  \url{https://doi.org/10.1090/ert/655}.

\bibitem{ah2}
{\sc P.~N. Achar and W.~Hardesty}, {\em Silting complexes of coherent sheaves
  and the {H}umphreys conjecture}, Duke Math. J., 173 (2024), pp.~2397--2445,
  \url{https://doi.org/10.1215/00127094-2023-0060}.

\bibitem{ahr-conj}
{\sc P.~N. Achar, W.~Hardesty, and S.~Riche}, {\em Conjectures on tilting
  modules and antispherical $p$-cells}.
\newblock Preprint~\href{https://arxiv.org/abs/1812.09960}{arXiv:1812.09960}.

\bibitem{ahr}
{\sc P.~N. Achar, W.~Hardesty, and S.~Riche}, {\em On the {H}umphreys
  conjecture on support varieties of tilting modules}, Transform. Groups, 24
  (2019), pp.~597--657, \url{https://doi.org/10.1007/s00031-019-09513-y}.

\bibitem{ahr2}
{\sc P.~N. Achar, W.~Hardesty, and S.~Riche}, {\em Representation theory of
  disconnected reductive groups}, Doc. Math., 25 (2020), pp.~2149--2177,
  \url{https://doi.org/10.4171/dm/796}.

\bibitem{ahr3}
{\sc P.~N. Achar, W.~Hardesty, and S.~Riche}, {\em Integral exotic sheaves and
  the modular {L}usztig--{V}ogan bijection}, J. Lond. Math. Soc. (2), 106
  (2022), pp.~2403--2458, \url{https://doi.org/10.1112/jlms.12638}.

\bibitem{amrw}
{\sc P.~N. Achar, S.~Makisumi, S.~Riche, and G.~Williamson}, {\em Koszul
  duality for {K}ac--{M}oody groups and characters of tilting modules}, J.
  Amer. Math. Soc., 32 (2019), pp.~261--310,
  \url{https://doi.org/10.1090/jams/905}.

\bibitem{ar-reductive}
{\sc P.~N. Achar and S.~Riche}, {\em Reductive groups, the loop {G}rassmannian,
  and the {S}pringer resolution}, Invent. Math., 214 (2018), pp.~289--436,
  \url{https://doi.org/10.1007/s00222-018-0805-1}.

\bibitem{ar-steinberg}
{\sc P.~N. Achar and S.~Riche}, {\em A geometric {S}teinberg formula},
  Transform. Groups, 28 (2023), pp.~1001--1032,
  \url{https://doi.org/10.1007/s00031-022-09768-y}.

\bibitem{ar-blocks}
{\sc P.~N. Achar and S.~Riche}, {\em A geometric model for blocks of
  {F}robenius kernels}, Ark. Mat., 62 (2024), pp.~217--329,
  \url{https://doi.org/10.4310/arkiv.2024.v62.n2.a1}.

\bibitem{andersen}
{\sc H.~H. Andersen}, {\em The strong linkage principle}, J. Reine Angew.
  Math., 315 (1980), pp.~53--59,
  \url{https://doi.org/10.1515/crll.1980.315.53}.

\bibitem{andersen-tilting}
{\sc H.~H. Andersen}, {\em Tilting modules for algebraic groups}, in Algebraic
  groups and their representations ({C}ambridge, 1997), vol.~517 of NATO Adv.
  Sci. Inst. Ser. C: Math. Phys. Sci., Kluwer Acad. Publ., Dordrecht, 1998,
  pp.~25--42.

\bibitem{andersen-jantzen}
{\sc H.~H. Andersen and J.~C. Jantzen}, {\em Cohomology of induced
  representations for algebraic groups}, Math. Ann., 269 (1984), pp.~487--525,
  \url{https://doi.org/10.1007/BF01450762}.

\bibitem{bnps2}
{\sc C.~P. Bendel, D.~K. Nakano, C.~Pillen, and P.~Sobaje}, {\em On {D}onkin's
  tilting module conjecture {II}: counterexamples}, Compos. Math., 160 (2024),
  pp.~1167--1193, \url{https://doi.org/10.1112/S0010437X24007115}.

\bibitem{bnps3}
{\sc C.~P. Bendel, D.~K. Nakano, C.~Pillen, and P.~Sobaje}, {\em On {D}onkin's
  {T}ilting {M}odule {C}onjecture {III}: {N}ew generic lower bounds}, J.
  Algebra, 655 (2024), pp.~95--109,
  \url{https://doi.org/10.1016/j.jalgebra.2023.07.001}.

\bibitem{bez-qes}
{\sc R.~Bezrukavnikov}, {\em Quasi-exceptional sets and equivariant coherent
  sheaves on the nilpotent cone}, Represent. Theory, 7 (2003), pp.~1--18,
  \url{https://doi.org/10.1090/S1088-4165-03-00158-4}.

\bibitem{bez-tcac}
{\sc R.~Bezrukavnikov}, {\em On tensor categories attached to cells in affine
  {W}eyl groups}, in Representation theory of algebraic groups and quantum
  groups, vol.~40 of Adv. Stud. Pure Math., Math. Soc. Japan, Tokyo, 2004,
  pp.~69--90, \url{https://doi.org/10.2969/aspm/04010069}.

\bibitem{bezrukavnikov}
{\sc R.~Bezrukavnikov}, {\em Cohomology of tilting modules over quantum groups
  and {$t$}-structures on derived categories of coherent sheaves}, Invent.
  Math., 166 (2006), pp.~327--357,
  \url{https://doi.org/10.1007/s00222-006-0514-z}.

\bibitem{bez-psaf}
{\sc R.~Bezrukavnikov}, {\em Perverse sheaves on affine flags and nilpotent
  cone of the {L}anglands dual group}, Israel J. Math., 170 (2009),
  pp.~185--206, \url{https://doi.org/10.1007/s11856-009-0025-x}.

\bibitem{br-2real}
{\sc R.~Bezrukavnikov and S.~Riche}, {\em On two modular geometric realizations
  of an affine {H}ecke algebra}.
\newblock Preprint~\href{https://arxiv.org/abs/2402.08281}{arXiv:2402.08281}.

\bibitem{br-Hecke}
{\sc R.~Bezrukavnikov and S.~Riche}, {\em Hecke action on the principal block},
  Compos. Math., 158 (2022), pp.~953--1019,
  \url{https://doi.org/10.1112/s0010437x22007436}.

\bibitem{ciappara}
{\sc J.~Ciappara}, {\em Hecke category actions via {S}mith--{T}reumann theory},
  Compos. Math., 159 (2023), pp.~2089--2124,
  \url{https://doi.org/10.1112/s0010437x23007340}.

\bibitem{fm}
{\sc M.~Finkelberg and I.~Mirkovi\'{c}}, {\em Semi-infinite flags. {I}. {C}ase
  of global curve {$\bold P^1$}}, in Differential topology,
  infinite-dimensional {L}ie algebras, and applications, vol.~194 of Amer.
  Math. Soc. Transl. Ser. 2, Amer. Math. Soc., Providence, RI, 1999,
  pp.~81--112, \url{https://doi.org/10.1090/trans2/194/05}.

\bibitem{friedlander-parshall}
{\sc E.~M. Friedlander and B.~J. Parshall}, {\em Cohomology of {L}ie algebras
  and algebraic groups}, Amer. J. Math., 108 (1986), pp.~235--253,
  \url{https://doi.org/10.2307/2374473}.

\bibitem{friedlander-suslin}
{\sc E.~M. Friedlander and A.~Suslin}, {\em Cohomology of finite group schemes
  over a field}, Invent. Math., 127 (1997), pp.~209--270,
  \url{https://doi.org/10.1007/s002220050119}.

\bibitem{ginzburg-kumar}
{\sc V.~Ginzburg and S.~Kumar}, {\em Cohomology of quantum groups at roots of
  unity}, Duke Math. J., 69 (1993), pp.~179--198,
  \url{https://doi.org/10.1215/S0012-7094-93-06909-8}.

\bibitem{hardesty}
{\sc W.~D. Hardesty}, {\em On support varieties and the {H}umphreys conjecture
  in type {$A$}}, Adv. Math., 329 (2018), pp.~392--421,
  \url{https://doi.org/10.1016/j.aim.2018.01.023}.

\bibitem{humphreys}
{\sc J.~E. Humphreys}, {\em Comparing modular representations of semisimple
  groups and their {L}ie algebras}, in Modular interfaces ({R}iverside, {CA},
  1995), vol.~4 of AMS/IP Stud. Adv. Math., Amer. Math. Soc., Providence, RI,
  1997, pp.~69--80, \url{https://doi.org/10.1090/amsip/004/05}.

\bibitem{jantzen}
{\sc J.~C. Jantzen}, {\em Representations of algebraic groups}, vol.~107 of
  Mathematical Surveys and Monographs, American Mathematical Society,
  Providence, RI, second~ed., 2003.

\bibitem{jantzen-char}
{\sc J.~C. Jantzen}, {\em Character formulae from {H}ermann {W}eyl to the
  present}, in Groups and analysis, vol.~354 of London Math. Soc. Lecture Note
  Ser., Cambridge Univ. Press, Cambridge, 2008, pp.~232--270,
  \url{https://doi.org/10.1017/CBO9780511721410.012}.

\bibitem{jensen-williamson}
{\sc L.~T. Jensen and G.~Williamson}, {\em The {$p$}-canonical basis for
  {H}ecke algebras}, in Categorification and higher representation theory,
  vol.~683 of Contemp. Math., Amer. Math. Soc., Providence, RI, 2017,
  pp.~333--361, \url{https://doi.org/10.1090/conm/683}.

\bibitem{kazhdan-lusztig}
{\sc D.~Kazhdan and G.~Lusztig}, {\em Representations of {C}oxeter groups and
  {H}ecke algebras}, Invent. Math., 53 (1979), pp.~165--184,
  \url{https://doi.org/10.1007/BF01390031}.

\bibitem{lusztig-pbs}
{\sc G.~Lusztig}, {\em Some problems in the representation theory of finite
  {C}hevalley groups}, in The {S}anta {C}ruz {C}onference on {F}inite {G}roups
  ({U}niv. {C}alifornia, {S}anta {C}ruz, {C}alif., 1979), vol.~37 of Proc.
  Sympos. Pure Math., Amer. Math. Soc., Providence, RI, 1980, pp.~313--317.

\bibitem{lusztig-cawg4}
{\sc G.~Lusztig}, {\em Cells in affine {W}eyl groups. {IV}}, J. Fac. Sci. Univ.
  Tokyo Sect. IA Math., 36 (1989), pp.~297--328.

\bibitem{mv}
{\sc I.~Mirkovi\'{c} and K.~Vilonen}, {\em Geometric {L}anglands duality and
  representations of algebraic groups over commutative rings}, Ann. of Math.
  (2), 166 (2007), pp.~95--143,
  \url{https://doi.org/10.4007/annals.2007.166.95}.

\bibitem{npv}
{\sc D.~K. Nakano, B.~J. Parshall, and D.~C. Vella}, {\em Support varieties for
  algebraic groups}, J. Reine Angew. Math., 547 (2002), pp.~15--49,
  \url{https://doi.org/10.1515/crll.2002.049}.

\bibitem{ostrik-tensor}
{\sc V.~Ostrik}, {\em Tensor ideals in the category of tilting modules},
  Transform. Groups, 2 (1997), pp.~279--287,
  \url{https://doi.org/10.1007/BF01234661}.

\bibitem{ostrik}
{\sc V.~V. Ostrik}, {\em Cohomological supports for quantum groups},
  Funktsional. Anal. i Prilozhen., 32 (1998), pp.~22--34, 95,
  \url{https://doi.org/10.1007/BF02463206}.

\bibitem{riche}
{\sc S.~Riche}, {\em Lectures on modular representation theory of reductive
  algebraic groups}.
\newblock Book in preparation, preliminary version available on the author's
  webpage.

\bibitem{riche-survey}
{\sc S.~Riche}, {\em Some applications of the geometric {S}atake equivalence to
  modular representation theory}.
\newblock Preprint~\href{https://arxiv.org/abs/2403.03734}{arXiv:2403.03734}.

\bibitem{rw-tilting}
{\sc S.~Riche and G.~Williamson}, {\em Tilting modules and the {$p$}-canonical
  basis}, Ast\'{e}risque,  (2018), pp.~ix+184,
  \url{https://doi.org/10.24033/ast.1043}.

\bibitem{rw-simple}
{\sc S.~Riche and G.~Williamson}, {\em A simple character formula}, Ann. H.
  Lebesgue, 4 (2021), pp.~503--535, \url{https://doi.org/10.5802/ahl.79}.

\bibitem{rw-smith}
{\sc S.~Riche and G.~Williamson}, {\em Smith--{T}reumann theory and the linkage
  principle}, Publ. Math. Inst. Hautes \'{E}tudes Sci., 136 (2022),
  pp.~225--292, \url{https://doi.org/10.1007/s10240-022-00134-y}.

\bibitem{sobaje}
{\sc P.~Sobaje}, {\em On character formulas for simple and tilting modules},
  Adv. Math., 369 (2020), pp.~107172, 8,
  \url{https://doi.org/10.1016/j.aim.2020.107172}.

\bibitem{soergel}
{\sc W.~Soergel}, {\em Kazhdan-{L}usztig polynomials and a combinatoric[s] for
  tilting modules}, Represent. Theory, 1 (1997), pp.~83--114,
  \url{https://doi.org/10.1090/S1088-4165-97-00021-6}.

\bibitem{verma}
{\sc D.-N. Verma}, {\em The r\^{o}le of affine {W}eyl groups in the
  representation theory of algebraic {C}hevalley groups and their {L}ie
  algebras}, in Lie groups and their representations ({P}roc. {S}ummer
  {S}chool, {B}olyai {J}\'{a}nos {M}ath. {S}oc., {B}udapest, 1971), Halsted
  Press, New York-Toronto, Ont., 1975, pp.~653--705.

\bibitem{vogan}
{\sc D.~A. Vogan, Jr.}, {\em The method of coadjoint orbits for real reductive
  groups}, in Representation theory of {L}ie groups ({P}ark {C}ity, {UT},
  1998), vol.~8 of IAS/Park City Math. Ser., Amer. Math. Soc., Providence, RI,
  2000, pp.~179--238, \url{https://doi.org/10.1090/pcms/008/05}.

\bibitem{williamson-alg}
{\sc G.~Williamson}, {\em Algebraic representations and constructible sheaves},
  Jpn. J. Math., 12 (2017), pp.~211--259,
  \url{https://doi.org/10.1007/s11537-017-1646-1}.

\bibitem{williamson-torsion}
{\sc G.~Williamson}, {\em Schubert calculus and torsion explosion}, J. Amer.
  Math. Soc., 30 (2017), pp.~1023--1046,
  \url{https://doi.org/10.1090/jams/868}.
\newblock With a joint appendix with Alex Kontorovich and Peter J. McNamara.

\bibitem{williamson}
{\sc G.~Williamson}, {\em Parity sheaves and the {H}ecke category}, in
  Proceedings of the {I}nternational {C}ongress of {M}athematicians---{R}io de
  {J}aneiro 2018. {V}ol. {I}. {P}lenary lectures, World Sci. Publ., Hackensack,
  NJ, 2018, pp.~979--1015.

\end{thebibliography}

\end{document}